\documentclass{article}
\usepackage{amsmath, amsthm, amssymb, tikz, hyperref,enumitem} 
\usepackage{xurl}
\usetikzlibrary{arrows.meta}
\usepackage[vcentermath,enableskew]{youngtab}
\usepackage{fullpage,verbatim}

\newtheorem{theorem}{Theorem}[section]
\newtheorem{proposition}[theorem]{Proposition}
\newtheorem{lemma}[theorem]{Lemma}
\newtheorem{corollary}[theorem]{Corollary}

\newtheorem{problem}[theorem]{Problem}

\theoremstyle{definition}
\newtheorem{definition}[theorem]{Definition}

\newtheorem{remark}[theorem]{Remark}

\renewcommand{\S}{\mathcal{S}}
\renewcommand{\P}{\mathcal{P}}
\newcommand{\F}{\mathcal{F}}

\newcommand{\oldlambda}{\psi}
\renewcommand{\L}{\mathcal{L}}
\newcommand{\T}{\mathcal{T}}
\newcommand{\A}{\mathbf{A}}
\DeclareMathOperator\des{des}
\DeclareMathOperator\Des{Des}
\DeclareMathOperator\asc{asc}
\DeclareMathOperator\Asc{Asc}
\DeclareMathOperator\hdes{hdes}
\DeclareMathOperator\HDes{HDes}

\DeclareMathOperator\SYT{SYT}
\DeclareMathOperator\Inv{Inv}

\DeclareMathOperator\rev{rev}
\DeclareMathOperator\flip{flip}
\DeclareMathOperator\comp{comp}
\newcommand\id{\operatorname{id}}

\DeclareMathOperator\lowlabels{\mathfrak{L}_{\downarrow}}
\DeclareMathOperator\uplabels{\mathfrak{L}^{\uparrow}}
\newcommand{\card}[1]{{\lvert #1 \rvert}}
\newcommand\AscC{\Asc^\Box}
\newcommand\HDesC{\HDes^\Box}
\newcommand{\La}{\Lambda}
\newcommand{\Ch}{\chi}
\newcommand\DesC{\Des^\Box}
\newcommand{\row}{\rho}
\renewcommand{\skew}{{\lambda/\mu}}
\newcommand{\one}{\textcolor{blue}{1}}
\newcommand{\two}{\textcolor{blue}{2}}
\newcommand{\three}{\textcolor{blue}{3}}
\newcommand{\four}{\textcolor{blue}{4}}
\newcommand{\five}{\textcolor{blue}{5}}
\newcommand{\six}{\textcolor{blue}{6}}
\newcommand{\seven}{\textcolor{blue}{7}}
\newcommand{\eight}{\textcolor{blue}{8}}
\newcommand{\nine}{\textcolor{blue}{9}}

\newcommand{\drawtab}[2]{
\foreach \ro [count=\i] in #2
	\foreach \num [count=\j] in \ro
		{\draw[gray] (\j-1,#1-\i) rectangle (\j,#1-\i+1);
		\node[scale=1,blue] at (\j-.5,#1+.5-\i) {$\num$};}
}

\def\e{-- ++(1,0)}
\def\n{-- ++(0,1)}
\def\di{-- ++(.5,.5)}
\def\shortdi{-- ++(.3,.3)}
\newcommand\drawcircpath[3]{
\draw[#1] #2 circle (.3) \di#3\di circle (.3);
\draw[#1,ultra thick] #2++(.2,.2) \shortdi#3\shortdi;
}
\newcommand\drawcircpathsmall[3]{
\draw[#1] #2 circle (.3) \di#3\di circle (.27);
\draw[#1,ultra thick] #2++(.2,.2) \shortdi#3\shortdi;
}
\newcommand\drawcircpathbig[3]{
\draw[#1] #2 circle (.33) \di#3\di circle (.3);
\draw[#1,ultra thick] #2++(.2,.2) \shortdi#3\shortdi;
}
\newcommand\drawpath[3]{
\draw[#1,ultra thick] #2 circle (.04) \di#3\di circle (.04);
}
\newcommand\drawpathsmall[3]{
\draw[#1,ultra thick] #2 circle (.02) \di#3\di circle (.04);
}

\newcommand{\drawskewtab}[3]{ 
\foreach \ro [count=\i] in #3
	\foreach \num [count=\j] in \ro
		{\draw[gray] (#2[\i-1]+\j-1,#1-\i) rectangle (#2[\i-1]+\j,#1-\i+1);
		\node[scale=1,blue] at (#2[\i-1]+\j-.5,#1+.5-\i) {$\num$};}
}
\newcommand{\drawskewtabblack}[3]{ 
\foreach \ro [count=\i] in #3
	\foreach \num [count=\j] in \ro
		{\draw[gray] (#2[\i-1]+\j-1,#1-\i) rectangle (#2[\i-1]+\j,#1-\i+1);
		\node[scale=.8,black] at (#2[\i-1]+\j-.8,#1+.8-\i) {$\num$};}
}

\newcommand{\drawferrersdiagram}[3]{ 
\foreach \ro [count=\i] in #3
	\foreach \num [count=\j] in \ro
		{\filldraw (#2[\i-1]+\j-.5,#1-\i+.5) circle (2pt);
		\node[scale=.9,black] at (#2[\i-1]+\j-.7,#1+.8-\i) {$\num$};}
}

\title{Rowmotion on standard Young tableaux via chutes and ladders}
\author{Sergi Elizalde~\thanks{Department of Mathematics, Dartmouth College, Hanover, NH 03755, USA. E-mail: {\tt sergi.elizalde@dartmouth.edu}. 
} 
\and Jessica Striker~\thanks{Department of Mathematics, North Dakota State University, Fargo, ND 58108, USA. E-mail: {\tt jessica.striker@ndsu.edu}. 
} 
}

\begin{document}

\maketitle



\begin{abstract}
We introduce a rowmotion operation on standard Young tableaux of any given (straight or skew) shape. This operation generalizes the rowmotion map on  standard Young tableaux of shape $2\times m$ obtained by interpreting them as order ideals of the Type $A_{m-1}$ root poset, and the more general rowmotion map on $k\times m$ standard Young tableaux recently proposed by Elizalde using so-called arrow encodings.

We give two  definitions of our rowmotion operation, which we prove are equivalent. The first description is a direct bijection on standard Young tableaux that takes certain ascending lattice paths inside the tableaux (which we call ladders) to certain descending paths (which we call chutes); in particular, it sends the number of ascents to the number of high descents. Our lattice paths generalize Reading's arc diagrams for permutations.
The second description defines a partial order on standard Young tableaux of a given shape, and then uses Barnard's notion of rowmotion on semidistributive lattices. 

We use this correspondence to prove an alternative description of rowmotion on the weak order. We also study for which shapes rowmotion and promotion on standard Young tableaux have the same orbit structure.
\end{abstract}

\section{Introduction}
\label{sec:intro}
Rowmotion is a fundamental operation in dynamical algebraic combinatorics. Originally defined on finite distributive lattices, or equivalently, order ideals of posets \cite{brouwer1974period,Cfdf,Panyushev2009,striker2017dac,prorow}, 
rowmotion has recently been extended to lattices with weaker distributive-type properties \cite{Barnard2019,ThomasWilliams2019,Semidistrim} and objects that generalize order ideals \cite{ICS,ICS2row}. 

The goal of this paper is to define a rowmotion operation on the set $\SYT(\skew)$ of standard Young tableaux of shape $\skew$, for any given skew shape. The operation that we define is natural in the sense that, as we will show, one arrives at the same map in multiple ways. Additionally, when restricted to standard Young tableaux of rectangular shape $2\times m$, our operation coincides with the rowmotion map that arises when interpreting these as Dyck paths, or equivalently, as order ideals of the Type $A_{m-1}$ root poset; see Figure~\ref{fig:2row}.

The first way to describe this operation is by directly defining a bijection on standard Young tableaux that behaves well with respect to certain statistics; specifically, it takes the number of ascents to the number of high descents. By adding some additional information to these statistics --- keeping track of what we call {\em chutes} and {\em ladders} --- we show that there is a unique bijection with the desired property (see Theorem~\ref{thm:La} and Definition~\ref{def:row}).

A second way to arrive at our rowmotion operation starts by defining a partial order on $\SYT(\skew)$, whose cover relations are given by Bender--Knuth involutions, i.e., swaps of entries with consecutive values. This poset is a lattice that may be interpreted as a lower interval in weak Bruhat order, determined by the shape $\skew$. In the case where this shape is a $2\times m$ rectangle, this poset is the distributive lattice of order ideals of the Type $A_{m-1}$ root poset, but in general it is merely a \emph{semidistributive} lattice (see Corollary~\ref{cor:semidistr} and Theorem~\ref{thm:dist}), so it cannot be the lattice of order ideals of some other poset. 
This lattice appears in work of Bj\"orner and Wachs~\cite{Bjorner_Wachs}, where it is an example of a generalized quotient of a Coxeter group, and in more recent work of Lacina \cite[Sec.~6]{Lacina}, where it arises as a \emph{maximal chain descent order} of an interval in Young’s lattice. 

A notion of rowmotion on semidistributive lattices was developed by Barnard~\cite{Barnard2019}, with subsequent extensions to the related families of trim lattices by Thomas and Williams~\cite{ThomasWilliams2019} and semidistrim lattices by Defant and Williams~\cite{Semidistrim}. 
We also use ideas of Reading \cite{Reading2015}, which studies the canonical join complex of the weak order on permutations by representing them as \emph{arc diagrams}.
We show in Theorem~\ref{thm:kappa} that, perhaps surprisingly, the map resulting from this lattice-theoretic approach coincides with the rowmotion operation that we had defined directly on the tableaux by using chutes and ladders. As a consequence of this result, we prove in Corollary~\ref{cor:rubey} that rowmotion on the weak order can alternatively be computed by reflecting the corresponding arc diagram and reversing the resulting permutation.

For standard Young tableaux of rectangular shape $2\times m$ and a few other shapes, the orbit structure of our rowmotion operation matches that of \emph{promotion}, the well-studied action on tableaux defined as a composition of Bender--Knuth involutions (see Corollary \ref{cor:match}).  But this is not true for general shapes. In particular, the order of promotion on tableaux of rectangular or staircase shape with $n$ cells divides $2n$, while for rowmotion, the order is generally much larger.  

Finally, for the case of standard Young tableaux of rectangular shape, our rowmotion operation coincides with the bijection described in~\cite[Prop.~5.3]{EliSYT} in terms of so-called {\em arrow encodings} of tableaux, which in turn can be viewed as a generalization of rowmotion on Dyck paths.

The paper is structured as follows. 
In Section~\ref{sec:background}, we give some background on standard Young tableaux, and discuss the cases in which a rowmotion operation is known. We also review the noncrossing arc diagrams introduced by Reading~\cite{Reading2015}.
In Section~\ref{sec:rowmotion-chutes&ladders}, we describe two encodings of standard Young tableaux as families of noncrossing lattice paths. 
These paths, which we call {\em chutes} and {\em ladders}, refine high descents and ascents on standard Young tableaux, respectively, and they 
can be viewed as a generalization of Reading's noncrossing arc diagrams.
We then use the chutes and ladders to define a rowmotion operation on standard Young tableaux of any given skew shape. 
In Section~\ref{sec:kappa} we describe a semidistributive lattice on standard Young tableaux of a given shape, which is a sublattice of the weak order on permutations. We prove our main theorem showing that our rowmotion map on standard Young tableaux is equivalent to rowmotion on this semidistributive lattice; we then derive some consequences.
In Section~\ref{sec:open}, we state some open problems and possible directions for future research.

\section{Background on standard Young tableaux, rowmotion, and noncrossing arc diagrams}
\label{sec:background}
\subsection{Standard Young tableaux}\label{sec:SYT}

Fix a partition $\lambda=(\lambda_1,\dots,\lambda_k)$, where $\lambda_1\ge\lambda_2\ge\dots\ge\lambda_k>0$. We will draw its Young diagram in English notation, so the $i$th row from the top has $\lambda_i$ cells.
Fix another (possibly empty) partition $\mu$ whose Young diagram is contained in that of $\lambda$, and write it as
$\mu=(\mu_1,\dots,\mu_k)$ by appending zeros if needed, where $\mu_1\ge\mu_2\ge\dots\ge\mu_k\ge0$. We have $\mu_r\le \lambda_r$ for $1\le r\le k$.
Denote by $\skew$  the corresponding skew shape, and let $n=\sum_{r=1}^k (\lambda_r-\mu_r)$ be its number of cells. When $\mu$ is empty, we simply write $\lambda$ instead of $\lambda/\mu$.

For $1\le r\le k$ and $\mu_r< c\le\lambda_r$, we write $(r,c)$ to denote the cell in row $r$ and column $c$.
Denote by $\SYT(\skew)$ the set of standard Young tableaux of shape $\skew$. These are fillings of the cells of the Young diagram of $\skew$ with the numbers $1,2,\dots,n$ so that each number appears exactly once, the entries in rows increase from left to right, and the entries in columns increase from top to bottom. A special case of interest is when $\mu$ is the empty partition, in which case $\SYT(\lambda)$ denotes the set of standard Young tableaux of straight shape $\lambda$.

\begin{definition}\label{def:ascents&descents}
Let $T\in\SYT(\skew)$. We will denote by $C_i=(r_i,c_i)$ the cell in $T$ that contains $i$. For $1\le i\le n-1$, we say that $i$ is
\begin{itemize}
\item an {\em ascent} if $r_i>r_{i+1}$ (i.e., $i+1$ is higher than $i$),
\item a {\em descent} if $r_i<r_{i+1}$ (i.e., $i+1$ is lower than $i$),
\item a {\em high descent} if $r_i<r_{i+1}$ and $c_i>c_{i+1}$ (i.e., $i+1$ is lower than and to the left of $i$).
\end{itemize}
Denote by $\Asc(T)$, $\Des(T)$ and $\HDes(T)$ the sets of ascents, descents and high descents of $T$, respectively. 
Let $\asc(T)=\card{\Asc(T)}$ and $\hdes(T)=\card{\HDes(T)}$.
\end{definition}

Let $T_R$ denote the unique tableau in $\SYT(\skew)$ where all the entries in each row are smaller than all the entries in the next row. 
Specifically, $T_R$ has entries $1,2,\dots,\lambda_1-\mu_1$ in the first row, entries $\lambda_1-\mu_1+1,\lambda_1-\mu_1+2,\dots,\lambda_1-\mu_1+\lambda_2-\mu_2$ in the second row, and so on. 
Similarly, let $T_C$ denote the unique tableau in $\SYT(\skew)$ where all the entries in each column are smaller than all the entries in the next column. 
For example, for shape $\lambda=(4,3,2)$, we have
$$T_R=\young(1234,567,89),\quad T_C=\young(1479,258,36),$$
and for shape $\lambda/\mu=(4,3,2)/(1,1,0)$, we have
$$T_R=\young(:123,:45,67),\quad T_C=\young(:257,:36,14).$$

\begin{definition}\label{def:labeling_by_rows}
    The {\em labeling by rows} (resp., {\em by columns}) 
    of a shape $\skew$ assigns to each cell the label given by its entry in $T_R$ (resp., $T_C$). 
\end{definition}

\subsection{Rowmotion on $\SYT(m^2)$ and the Type $A_{m-1}$ positive root poset} 
\label{sec:two_row}
It is well-known that $\card{\SYT(m^2)}=\frac{1}{m+1}\binom{2m}{m}$, the $m$th Catalan number. 
A well-known bijection between the set $\SYT(m^2)$ and the set of Dyck paths of semilength $m$, i.e., lattice paths from $(0,0)$ to $(2m,0)$ comprised of $m$ up-steps $(1,1)$ and $m$ down-steps $(1,-1)$ that stay weakly above the $x$-axis, proceeds as follows.
Given $T\in\SYT(m^2)$, define  the $i$th step in the resulting Dyck path $D$ to be an up-step if $i$ is in the top row of $T$, and a down-step if $i$ is in the bottom row of $T$. The condition that the Dyck path stays weakly above the $x$-axis translates to the column-strict condition on the standard Young tableau. 

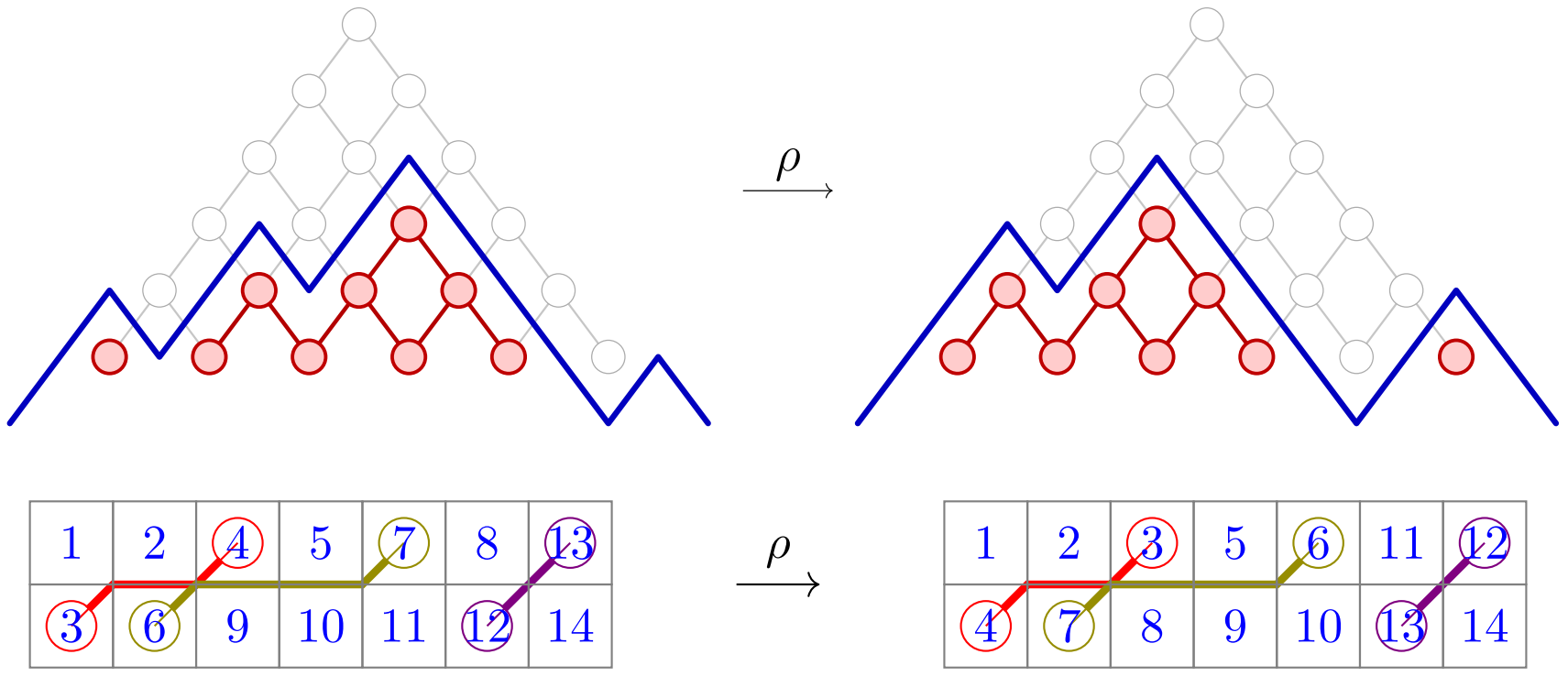
\begin{figure}[htbp]
    \centering
    \scalebox{.6}{
\begin{tikzpicture}[
    x=1.2cm,
    y=0.8cm,
    edge/.style={
        draw=gray!45,
        line width=0.7pt
    },
    ideal edge/.style={
        draw=red!70!black,
        line width=1.4pt
    },
    root/.style={
        circle,
        draw=gray!65,
        fill=white,
        minimum size=4mm,
        inner sep=0pt
    },
    ideal/.style={
        root,
        draw=red!75!black,
        fill=red!20,
        line width=1.2pt
    },
    dyck path/.style={
        draw=blue!75!black,
        line width=2pt,
        line join=round,
        line cap=round
    },
    row arrow/.style={
        -{Latex[length=3mm,width=2mm]},
        line width=1.2pt
    }
]

\coordinate (r11) at (0,0);
\coordinate (r22) at (1,0);
\coordinate (r33) at (2,0);
\coordinate (r44) at (3,0);
\coordinate (r55) at (4,0);
\coordinate (r66) at (5,0);

\coordinate (r12) at (0.5,1);
\coordinate (r23) at (1.5,1);
\coordinate (r34) at (2.5,1);
\coordinate (r45) at (3.5,1);
\coordinate (r56) at (4.5,1);

\coordinate (r13) at (1,2);
\coordinate (r24) at (2,2);
\coordinate (r35) at (3,2);
\coordinate (r46) at (4,2);

\coordinate (r14) at (1.5,3);
\coordinate (r25) at (2.5,3);
\coordinate (r36) at (3.5,3);

\coordinate (r15) at (2,4);
\coordinate (r26) at (3,4);

\coordinate (r16) at (2.5,5);

\foreach \a/\b in {r11/r12,r22/r12,r22/r23,r33/r23,r33/r34,r44/r34,r44/r45,r55/r45,r55/r56,r66/r56,r12/r13,r23/r13,r23/r24,r34/r24,r34/r35,r45/r35,r45/r46,r56/r46,r13/r14,r24/r14,r24/r25,r35/r25,r35/r36,r46/r36,r14/r15,r25/r15,r25/r26,r36/r26,r15/r16,r26/r16}{
    \draw[edge] (\a) -- (\b);
}

\foreach \a/\b in {r22/r23,r33/r23,r33/r34,r44/r34,r44/r45,r55/r45,r34/r35,r45/r35}{
    \draw[ideal edge] (\a) -- (\b);
}

\draw[dyck path]
    (-1,-1)
    -- (-0.5,0)
    -- (0,1)
    -- (0.5,0)
    -- (1,1)
    -- (1.5,2)
    -- (2,1)
    -- (2.5,2)
    -- (3,3)
    -- (3.5,2)
    -- (4,1)
    -- (4.5,0)
    -- (5,-1)
    -- (5.5,0)
    -- (6,-1);

\foreach \v in {r11,r22,r33,r44,r55,r66,r12,r23,r34,r45,r56,r13,r24,r35,r46,r14,r25,r36,r15,r26,r16}{
    \node[root] at (\v) {};
}

\foreach \v in {r11,r22,r33,r44,r55,r23,r34,r45,r35}{
    \node[ideal] at (\v) {};
}

\draw[->] (6.35,2.5) -- (7.25,2.5)
    node[midway,above] {\LARGE$\rho$};

\coordinate (s11) at (8.5,0);
\coordinate (s22) at (9.5,0);
\coordinate (s33) at (10.5,0);
\coordinate (s44) at (11.5,0);
\coordinate (s55) at (12.5,0);
\coordinate (s66) at (13.5,0);

\coordinate (s12) at (9,1);
\coordinate (s23) at (10,1);
\coordinate (s34) at (11,1);
\coordinate (s45) at (12,1);
\coordinate (s56) at (13,1);

\coordinate (s13) at (9.5,2);
\coordinate (s24) at (10.5,2);
\coordinate (s35) at (11.5,2);
\coordinate (s46) at (12.5,2);

\coordinate (s14) at (10,3);
\coordinate (s25) at (11,3);
\coordinate (s36) at (12,3);

\coordinate (s15) at (10.5,4);
\coordinate (s26) at (11.5,4);

\coordinate (s16) at (11,5);

\foreach \a/\b in {s11/s12,s22/s12,s22/s23,s33/s23,s33/s34,s44/s34,s44/s45,s55/s45,s55/s56,s66/s56,s12/s13,s23/s13,s23/s24,s34/s24,s34/s35,s45/s35,s45/s46,s56/s46,s13/s14,s24/s14,s24/s25,s35/s25,s35/s36,s46/s36,s14/s15,s25/s15,s25/s26,s36/s26,s15/s16,s26/s16}{
    \draw[edge] (\a) -- (\b);
}

\foreach \a/\b in {s11/s12,s22/s12,s22/s23,s33/s23,s33/s34,s44/s34,s23/s24,s34/s24}{
    \draw[ideal edge] (\a) -- (\b);
}

\draw[dyck path]
    (7.5,-1)
    -- (8,0)
    -- (8.5,1)
    -- (9,2)
    -- (9.5,1)
    -- (10,2)
    -- (10.5,3)
    -- (11,2)
    -- (11.5,1)
    -- (12,0)
    -- (12.5,-1)
    -- (13,0)
    -- (13.5,1)
    -- (14,0)
    -- (14.5,-1);

\foreach \v in {s11,s22,s33,s44,s55,s66,s12,s23,s34,s45,s56,s13,s24,s35,s46,s14,s25,s36,s15,s26,s16}{
    \node[root] at (\v) {};
}

\foreach \v in {s11,s22,s33,s44,s66,s12,s23,s34,s24}{
    \node[ideal] at (\v) {};
}

\end{tikzpicture}}

\vspace{.5cm}
\begin{tikzpicture}[scale=.6]
\drawcircpath{red}{(.5,.5)}{\e}
\drawcircpath{olive}{(1.5,.5)}{\e\e}
\drawcircpath{violet}{(5.5,.5)}{}
\drawtab{2}{{{1,2,4,5,7,8,13},{3,6,9,10,11,12,14}}}
\draw[->] (8.5,1) -- node[above]{$\rho$} (9.5,1);
\begin{scope}[shift={(11,0)}]
\drawcircpath{red}{(.5,.5)}{\e}
\drawcircpath{olive}{(1.5,.5)}{\e\e}
\drawcircpath{violet}{(5.5,.5)}{}
\drawtab{2}{{{1,2,3,5,6,11,12},{4,7,8,9,10,13,14}}}
\end{scope}
\end{tikzpicture}
    \caption{An example of rowmotion on an order ideal of the poset $\A_6$ and our rowmotion map $\row$, as defined in Section~\ref{sec:rowmotion}, 
    on the corresponding tableau of shape $\lambda=(7,7)$.}
    \label{fig:2row}
\end{figure}

Dyck paths may be interpreted as order ideals, i.e., downward-closed sets, of the Type $A_{m-1}$ positive root poset, which we denote by $\A_{m-1}$. This is the triangular poset with elements $\{(i,j) \ | \ 1\leq i<j\leq m\} $ and relations $(i,j)\leq(i',j')$ whenever $i\geq i'$ and $j\leq j'$. A straightforward bijection between order ideals of $\A_{m-1}$ and Dyck paths of semilength $m$ proceeds by drawing a lattice path that separates the order ideal from the rest of the poset; see for instance \cite[Sec.~1.2]{Fishel19} and Figure~\ref{fig:2row}.

Given a finite poset $P$, denote by $J(P)$ its set of order ideals.
Rowmotion on order ideals \cite{brouwer1974period,Panyushev2009,striker2017dac,prorow} is the bijection $\row:J(P)\rightarrow J(P)$ where, for any $I\in J(P)$, its image
$\row(I)$ is the order ideal generated by the minimal elements of $P\setminus I$. 
Using the above interpretation of order ideals of $\A_{m-1}$ as Dyck paths, and the above bijection between Dyck paths and $2$-row standard Young tableaux, one can think of rowmotion  on order ideals of  $\A_{m-1}$ as an operation on $\SYT(m^2)$.
We will see that 
the rowmotion action on $\SYT(\skew)$ that we define in Section~\ref{sec:rowmotion}, when restricted to the case $\skew=(m^2)$, coincides with this operation, justifying our use of the term rowmotion for this more general action. See Figure~\ref{fig:2row} for an example and Remark~\ref{remark:tworowrho} for further discussion.

Cameron and Fon-der-Flaass showed that rowmotion on order ideals can be equivalently defined in terms of \emph{toggles} \cite{Cfdf}. Given $x\in P$, the toggle $t_x:J(P)\rightarrow J(P)$ acting on $I$ returns the symmetric difference $I\triangle \{x\}$, provided this is an order ideal, and $I$ otherwise. Then $\row=t_{x_1}\circ t_{x_2}\circ\cdots\circ  t_{x_{|P|}}$, where the poset elements are naturally labeled in the order of a linear extension. That is, rowmotion toggles all elements of the poset from top to bottom.

As noted in \cite{prorow,stanley_proevac}, the operation on order ideals of $\A_{m-1}$  that toggles all the elements of this poset from left to right by \emph{column} (where the  $k$th column 
consists of all elements $(i,j)$ with $i+j=k+1$) corresponds to \emph{promotion} on the standard Young tableau associated to the order ideal. Promotion on a tableau in $\SYT(\skew)$   may be defined in terms of \emph{Bender--Knuth involutions}. The $k$th Bender--Knuth involution $b_k$ acts by swapping entries $k$ and $k+1$ if they are not in the same row and not in the same column, that is, if the result is still in $\SYT(\skew)$. Promotion is the composition of these involutions $b_{n-1}\circ\cdots\circ b_2\circ b_1$. The Bender--Knuth involution $b_k$ on $T\in\SYT(m^2)$ has the effect on the corresponding Dyck path of swapping steps $k$ and $k+1$, which is equivalent to toggling the elements in column $k$ of $\A_{m-1}$. Striker and Williams \cite{prorow} used the toggle group to show that toggling any poset from left to right by column (with a proper notion of column) has the same orbit structure as rowmotion. Thus, on $\SYT(m^2)$, promotion and rowmotion have the same orbit structure. We shall see that this correspondence between promotion and rowmotion orbit sizes does not extend to tableaux of general shape; see Corollary~\ref{cor:match} and Section~\ref{sec:orbits}. Nevertheless, we show in our main result, Theorem~\ref{thm:kappa}, that our notion of rowmotion on $\SYT(\skew)$ is natural in the sense that there is a semidistributive lattice whose $\kappa$-rowmotion operation (as defined in~\cite{Barnard2019}) coincides with our new tableau rowmotion~$\rho$.

\subsection{Rowmotion on $\SYT(m^k)$ via arrow encodings}
\label{sec:arrows}

A generalization of rowmotion on $\SYT(m^2)$ to standard Young tableaux of arbitrary rectangular shape  has been recently proposed in \cite[Prop.~5.3]{EliSYT}.
This operation is defined as a composition of certain maps that were constructed to bijectively prove the symmetry of the distribution of the number of descents on rectangular tableaux, as well as the relationship between the number of descents and the number of ascents.
These maps are described in terms of so-called {\em arrow encodings} of tableaux, which are configurations of up and down arrows associated to the cells of the Young diagram to indicate the relative location of each entry $i$ with respect to $i-1$.
We refer the reader to \cite{EliSYT} for details on these constructions.

It can be shown that the rowmotion map $\rho$ that we will introduce in Definition~\ref{def:row} in terms of chutes and ladders, when restricted to the set $\SYT(m^k)$ of rectangular tableaux, coincides with the rowmotion map proposed in \cite[Prop.~5.3]{EliSYT}.

\subsection{Noncrossing arc diagrams for permutations}\label{sec:arc_diagrams}

Let $\S_n$ denote the set of permutations of $\{1,2,\dots,n\}$. For $1\le i\le n-1$, we say that $i$ is a {\em descent} of $\pi\in\S_n$ if $\pi(i)>\pi(i+1)$, and that it is an {\em ascent} if $\pi(i)<\pi(i+1)$. 
Denote by $\Des(\pi)$ and $\Asc(\pi)$ the sets of descents and ascents of $\pi$, respectively. 

\begin{figure}[htbp]
    \centering
\begin{tikzpicture}[scale=.7,
    dot/.style={circle,fill=black,inner sep=2pt},
]

    \draw[red,thick] (0,2) to [bend left=55] (0,5);
    
    \draw[red,thick] (0,3)  to [bend right=55] (0,5.5) to [bend left=55] (0,8);

  \foreach \i in {1,...,8} {
    \node[dot] at (0,\i) {};
    \node[left=3pt] at (0,\i) {\(\i\)};
  }
\end{tikzpicture}
    \caption{The arc diagram of the permutation $15248367$.}
    \label{fig:arc_diagram}
\end{figure}
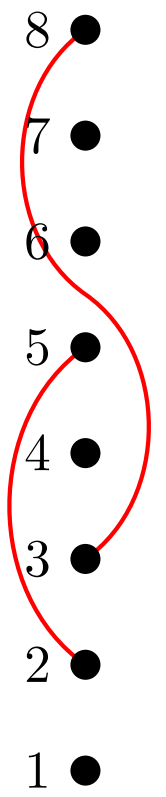

In \cite[Thm.~3.1]{Reading2015}, Reading describes a bijection between permutations and certain noncrossing arc diagrams.
The arc diagram of $\pi\in\S_n$ is constructed by first plotting a dot with label $k$ at coordinate $(0,k)$ for each $1\leq k\leq n$. Then, for each descent $b=\pi(i)>\pi(i+1)=a$ of $\pi$, 
draw an arc between dots $a$ and $b$ such that, for each $a<k<b$, the dot $k$ stays left or right of this arc according to whether $k$ appears to the left or to the right of $a$ and $b$ in $\pi$, respectively.
See Figure~\ref{fig:arc_diagram} for an example.

\section{Rowmotion on standard Young tableaux via chutes and ladders}
\label{sec:rowmotion-chutes&ladders}
In this section, we define two bijections between standard Young tableaux of a given shape $\skew$ and certain families of noncrossing paths. We then use these bijections to define a rowmotion map on standard Young tableaux. 

\subsection{Chutes and ladders}\label{sec:chutes&ladders}
\begin{figure}[htb]
\centering
\begin{tikzpicture}[scale=.8]
\drawcircpath{olive}{(1.5,2.5)}{}
\drawcircpath{red}{(.5,1.5)}{\e\e\n}
\drawcircpath{violet}{(.5,.5)}{\e}
\drawcircpath{orange}{(3.5,2.5)}{}
\drawcircpathsmall{green}{(2.5,.5)}{}
\drawcircpathbig{brown}{(3.5,1.5)}{\n\e}
\drawtab{4}{{{1,2,5,8,10,17},{3,4,6,9,18},{7,11,13,16},{12,14,15,19}}}
\draw (2.5,-.7) node {$\La(T)$};
\begin{scope}[shift={(8,0)}]
\drawcircpath{olive}{(.5,2.5)}{}
\drawcircpath{red}{(.5,1.5)}{\e}
\drawcircpathsmall{gray}{(.5,.5)}{}
\drawcircpathbig{violet}{(1.5,1.5)}{\e\e\n}
\drawcircpath{orange}{(1.5,.5)}{}
\drawcircpathsmall{green}{(3.5,.5)}{\n}
\drawcircpathbig{brown}{(4.5,2.5)}{}
\drawtab{4}{{{1,2,5,8,10,17},{3,4,6,9,18},{7,11,13,16},{12,14,15,19}}}
\draw (2.5,-.7) node {$\Ch(T)$};
\end{scope}
\end{tikzpicture}
\caption{The families of paths $\La(T)$ and $\Ch(T)$ for a tableau $T\in\SYT(6,5,4,4)$.}
\label{fig:La}
\end{figure}
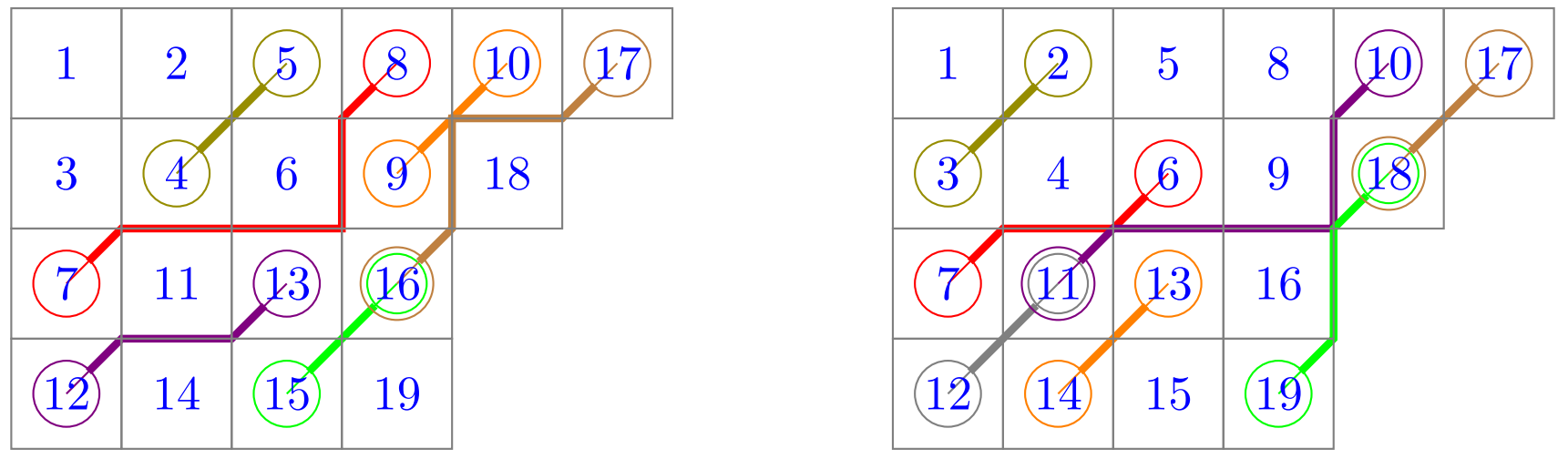

We start by defining two refinements of the statistics $\Asc$ and $\HDes$ from Section~\ref{sec:SYT}. The first refinement records the pairs of cells containing $i$ and $i+1$ for each ascent (resp.\ high descent) $i$ of $T\in\SYT(\skew)$:
\begin{align*} 
\AscC(T)&=\{\{C_i,C_{i+1}\}:i\in\Asc(T)\},\\
\HDesC(T)&=\{\{C_i,C_{i+1}\}:i\in\HDes(T)\}.
\end{align*}
For example, for the tableau $T$ in Figure~\ref{fig:La}, we have
\begin{align*} 
\AscC(T)&=\{\{(2,2),(1,3)\}, \{(3,1),(1,4)\}, \{(2,4),(1,5)\}, \{(4,1),(3,3)\}, \{(4,3),(3,4)\}, \{(3,4),(1,6)\}\},\\
\HDesC(T)&=\{\{(1,2),(2,1)\}, \{(2,3),(3,1)\}, \{(1,5),(3,2)\}, \{(3,2),(4,1)\}, \{(3,3),(4,2)\}, \{(1,6),(2,5)\},\\
&\hspace{6.5mm} \{(2,5),(4,4)\}\}.
\end{align*}

In general, the set of pairs $\AscC(T)$ does not uniquely determine the tableau $T\in\SYT(\skew)$. For example, the two tableaux  $$T=\young(\one\two\six,\three\four\seven,\five\eight\nine) \quad\text{and}\quad U=\young(\one\two\five,\three\six\seven,\four\eight\nine)$$
satisfy $\AscC(T)=\AscC(U)=\{\{(3,1),(1,3)\}\}$.

Let us introduce a further refinement of the statistic $\AscC$ that distinguishes these two tableaux.
Note that the entries strictly smaller than the ascent of $T$ form a tableau of shape $(2,2)$, whereas the entries strictly smaller than the ascent of $U$ form a tableau of shape $(2,1)$. Denote by $T^{<i}$ the standard Young tableau consisting of the cells of $T$ with entries strictly smaller than $i$, plus the cells in $\mu$. Our refined statistic will record the portion of the boundary of $T^{<i}$ between the cells where entries $i$ and $i+1$ lie in $T$, for each $i\in\Asc(T)$. As we will see later, this information is enough to recover $T$.

We will think of the cells of the Young diagram of $\skew$ as unit squares whose corners are points in $\mathbb{Z}^2$. 
Given two cells $C=(r,c)$ and $C'=(r',c')$, we say that $C'$ is north of (or higher than) $C$ if $r'<r$, and that $C'$ is west (or left) of $C$ if $c'<c$. We define south, lower, east, and right analogously.  We say that $C'$ is northeast of $C$ if both $r'<r$ and $c'>c$, and that $C'$ is weakly northeast of $C$ if $r'\le r$ and $c'\ge c$.

For pairs of cells where $C'$ is northeast of $C$, we will consider lattice paths from the northeast corner of $C$ to the southwest corner of $C'$ using unit steps north and east.
Note that such paths consist of $r-r'-1$ north steps and $c'-c-1$ east steps. 
We will slightly extend such paths at both endpoints, by adding a segment from the center of the cell $C$ to its northeast corner, and a segment from the southwest corner of the cell $C'$ to its center, so that each extended path starts at the center of $C$ and ends at the center of $C'$. Denote by $\P(C,C')$ the set of all such paths.
We define $\P(C',C)=\P(C,C')$, so that we can use this notation whenever one of the two cells is northeast of the other.

\begin{definition}\label{def:chutes&ladders}
For $i\in\Asc(T)$ or $i\in\HDes(T)$, let $P_i$ be the path in $\P(C_i,C_{i+1})$ whose north and east steps are determined by the boundary of $T^{<i}$. In other words, the cells above and to the left of the path are in $T^{<i}$, whereas those below and to the right of the path are not.
Define the following two sets of paths associated to $T$:
\[
\begin{array}{l@{\qquad}l}
\La(T)=\{P_i:i\in\Asc(T)\} & \text{(the {\em ladders} of $T$)},\\
\Ch(T)=\{P_i:i\in\HDes(T)\} & \text{(the {\em chutes} of $T$)}. 
\end{array}
\]
\end{definition}

An example is given in Figure~\ref{fig:La}. 
Note that $\card{\La(T)}=\asc(T)$ and $\card{\Ch(T)}=\hdes(T)$ by definition.

By construction, the paths in $\La(T)$ (resp.\ $\Ch(T)$) are noncrossing, even though they may intersect and even share steps. Additionally, each cell in $T$ can be the southwest cell of at most one path in the family, and the northeast cell of at most one path. Let us introduce some notation for families of paths with these properties.

\begin{definition}\label{def:families}
Let $\P_\skew=\bigcup_{C,C'}\P(C,C')$, where the union ranges over all pairs of cells $C,C'$ in $\skew$ with $C'$ northeast of $C$.
Let $\F_\skew$ be the set whose elements are families of paths $F\subseteq\P_\skew$ with the following three properties:
\begin{enumerate}[label=(\alph*)]
    \item No two paths in $F$ cross each other.
    \item No cell is the northeast cell of more than one path in $F$.
    \item No cell is the southwest cell of more than one path in $F$.
\end{enumerate}
\end{definition}

For any $T\in\SYT(\skew)$, we have $\La(T), \Ch(T)\in\F_\skew$ by construction.

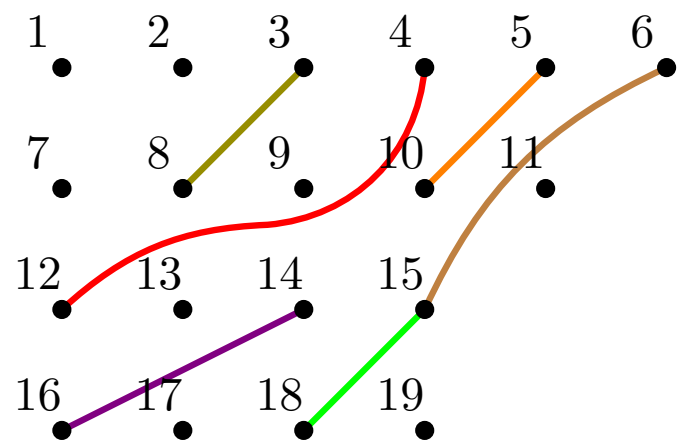
\begin{figure}[htb]
\centering
\begin{tikzpicture}[scale=.8]

\draw[olive,very thick] (1.5,2.5) to (2.5,3.5);
\draw[red,very thick] (0.5,1.5) to [bend left=20] (2.2,2.2) to [bend right=40] (3.5,3.5);
\draw[violet,very thick] (.5,.5) to (2.5,1.5);
\draw[orange,very thick] (3.5,2.5) to (4.5,3.5);
\draw[green,very thick] (2.5,.5) to (3.5,1.5);
\draw[brown,very thick] (3.5,1.5) to [bend left=20] (5.5,3.5);
\drawferrersdiagram{4}{{0,0,0,0}}{{{1,2,3,4,5,6},{7,8,9,10,11},{12,13,14,15},{16,17,18,19}}}
\end{tikzpicture}
\caption{The family of paths on the left of Figure~\ref{fig:La}, drawn as a noncrossing arc diagram.}
\label{fig:tableau_arc_diagram}
\end{figure}

\begin{remark}\label{rem:arc_diagrams}
There is an alternative way to think of the elements in $\F_\skew$, which is inspired by Reading's noncrossing arc diagrams for permutations discussed in Section~\ref{sec:arc_diagrams}. Instead of Young diagrams, consider Ferrers diagrams, where each cell is represented by the point at its center. A path in $\P(C,C')$, where $C'$ is northeast of $C$, can be viewed as a curve, called an {\em arc}, connecting the centers of the cells, moving monotonically in the northeast direction, and not passing through any other points of the Ferrers diagram. Two arcs in $\P(C,C')$ are equivalent if they pass to the right and to the left of the same points. We can then identify elements $F\in\F_\skew$ with families of arcs between points of the Ferrers diagram of $\skew$ satisfying the following three conditions:
\begin{enumerate}[label=(\alph*)]
    \item No two arcs in $F$ intersect, except possibly at their endpoints.
    \item No point is the northeast endpoint of more than one arc in $F$.
    \item No point is the southwest endpoint of more than one arc in $F$.
\end{enumerate}
We call such a family of arcs a {\em noncrossing arc diagram of shape $\skew$}. See Figure~\ref{fig:tableau_arc_diagram} for an example.
\end{remark}

Denote the staircase shape by $\delta_n=(n,n-1,\dots,1)$. 

\begin{remark}\label{rem:arc_diagrams-connection} In the special case that $\skew=\delta_n/\delta_{n-1}$, elements of $\F_{\delta_n/\delta_{n-1}}$ are equivalent to noncrossing arc diagrams for permutations from Section~\ref{sec:arc_diagrams}. 
Specifically, the cell in the $k$th row of the skew shape $\delta_n/\delta_{n-1}$ (which has label $k$ in the labeling by rows) corresponds to the point labeled $k$ in the arc diagram. Thus, a noncrossing arc diagram of shape $\delta_n/\delta_{n-1}$
can be turned into a noncrossing arc diagram in Reading's sense by  rotating it counterclockwise by 45$^\circ$ and then reflecting it vertically. For example, the arc diagram in Figure~\ref{fig:arc_diagram} corresponds to the tableau on the bottom left of Figure~\ref{fig:row_examples}.
\end{remark}

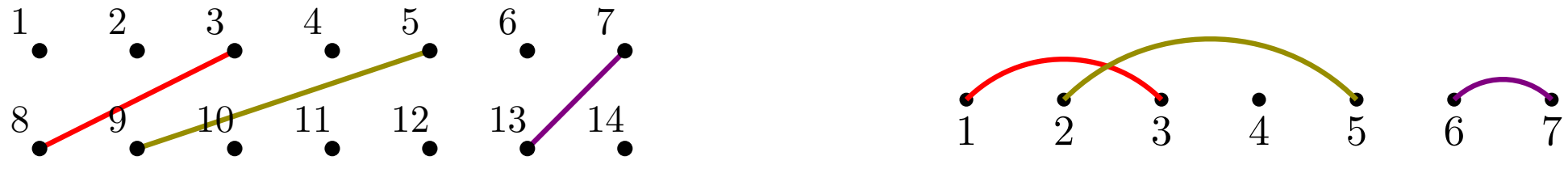
\begin{figure}[htb]
\centering
\begin{tikzpicture}[scale=.8]
\draw[red,very thick] (0.5,0.5) to (2.5,1.5);
\draw[olive,very thick] (1.5,.5) to (4.5,1.5);
\draw[violet,very thick] (5.5,.5) to (6.5,1.5);
\drawferrersdiagram{2}{{0,0,0,0}}{{{1,2,3,4,5,6,7},{8,9,10,11,12,13,14}}}
\begin{scope}[shift={(9,1)},every node/.style={draw,circle,fill,inner sep=1pt}]      
     \foreach \j in {1,2,3,4,5,6,7}
        \node[label=below:$\j$] at (\j,0)(\j) {};
\draw[red,very thick] (1,0) to [bend left=45] (3,0);
\draw[olive,very thick] (2,0) to [bend left=45] (5,0);
\draw[violet,very thick] (6,0) to [bend left=45] (7,0);
 \end{scope}
\end{tikzpicture}
\caption{The noncrossing arc diagram of shape $(7,7)$ from Figure~\ref{fig:2row}, and the corresponding partial nonnesting matching on $\{1,2,\dots,7\}$.}
\label{fig:nonnesting}
\end{figure}

\begin{remark}
    When $\skew=(m^2)$, there is a simple bijection between $\F_{(m^2)}$ and the set of partial nonnesting matchings on $\{1,2,\dots,m\}$. Given $F\in\F_{(m^2)}$, we obtain the corresponding matching as follows. For each path in $F$ with endpoints at cells $C=(2,c)$ and $C'=(1,c')$, where $1\le c<c'\le m$, we match $c$ and $c'$.
    See Figure~\ref{fig:nonnesting} for an example.
\end{remark}

The definitions of chutes and ladders are related to each other by conjugation, that is, reflection along the main diagonal. To make this precise, denote the conjugate of $\lambda$ by $\lambda'$, and the conjugate of $T\in\SYT(\skew)$ by $T'\in\SYT(\lambda'/\mu')$. There is a straightforward bijection between $\F_\skew$ and $\F_{\lambda'/\mu'}$ obtained by reflecting each path along the main diagonal. Denote the image of a family of paths $F\in\F_\skew$ by $F'\in\F_{\lambda'/\mu'}$. Then we can relate the definitions of chutes and ladders by
\begin{equation}\label{eq:Chutes-Ladders-transpose}
\Ch(T)'=\La(T').
\end{equation}

\subsection{From families of paths to standard Young tableaux}

\begin{figure}[htb]
\centering
\begin{tikzpicture}[scale=.8]
\drawskewtabblack{4}{{5,4,1,0}}{{{1,2},{3,4,5},{6,7,8,9,10},{11,12,13,14}}}
\draw[orange,ultra thick] (1.5,.5) circle (.04) \di\di circle (.04);
\draw[red,ultra thick] (2.5,1.5) circle (.02) \di\e\di circle (.04);
\draw[brown,ultra thick] (3.5,1.5) circle (.04)  \di\e\n\di circle (.04);
\draw[violet,ultra thick] (2.5,.5) circle (.04) \di\e\di circle (.04);
\draw[green,ultra thick] (5.5,1.5) circle (.04) \di\n\di circle (.04);
\foreach \x in {(1,1),(4,1),(2,2),(6,3),(7,3)}{
    \fill \x+(-.5,-.5) circle (.08);
}
\draw[brown] (4.9,3.4) node {$Q_1$};
\draw[orange!50!red] (2.7,2.2) node {$Q_2$};
\draw[black] (1.7,1.8) node {$Q_3$};
\end{tikzpicture}
\caption{A family of paths $F\in\F_\skew$, where $\lambda/\mu=(7,7,6,4)/(5,4,1)$, with the cells labeled by rows.}
\label{fig:F}
\end{figure}
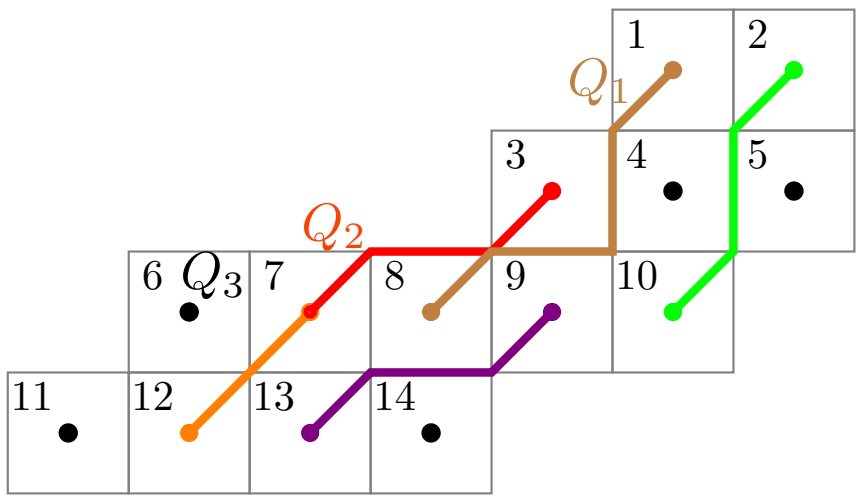

Fix a family of lattice paths $F\in\F_\skew$.
A {\em component} $Q$ of $F$ is a maximal set of cells $C^{(0)},C^{(1)},\dots,C^{(k)}$ (where $k\ge0$) and paths $P^{(1)},\dots,P^{(k)}\in F$ such that $P^{(j)}\in\P(C^{(j-1)},C^{(j)})$ for all $j$. Because of properties (b) and (c) in Definition~\ref{def:families}, every cell of $\skew$ belongs to exactly one component of $F$. By convention, $C^{(j)}$ is northeast of $C^{(j-1)}$ for all $j$, and we call $C^{(k)}$ the top cell of $Q$. In Figure~\ref{fig:F},
the component $Q_1$ consists of the brown path, and the component  $Q_2$ consists of the orange and red paths.

Given a component $Q$, define $L(Q)$ to be the set of cells of $\skew$ that are in the rows spanned by $Q$ and lie strictly to the left of $Q$. In Figure~\ref{fig:F}, $L(Q_1)$ consists of the cells labeled $3,6,7$, whereas $L(Q_2)$ consists of the cells labeled $6,11$. 
Similarly, define $A(Q)$ to be the set of cells of $\skew$ that are in the columns spanned by $Q$ and lie strictly above~$Q$. For example, $A(Q_2)$ is the cell labeled $6$.

\begin{figure}[htb]
\centering
\begin{tikzpicture}[scale=.6]
\begin{scope}[shift={(8,0)}]
\drawskewtab{5}{{3,2,1,0,0}}{{{3,9},{5,8},{4,6,11},{1,7},{2,10}}}
\end{scope}
\drawskewtab{5}{{3,2,1,0,0}}{{{,},{,},{,,},{,},{,}}}
\drawpath{brown}{(.5,.5)}{\n\n\e\n\e}
\drawpath{red}{(1.5,2.5)}{}
\drawpath{violet}{(1.5,1.5)}{-- ++(0,.05)\e--++(0,-.05)\n}
\drawpathsmall{purple}{(3.5,3.5)}{}
\drawpath{green}{(1.5,.5)}{-- ++(0,.95)\e-- ++(.0,.05)}
\foreach \x in {(1,2),(3,3)}{
    \fill \x+(-.5,-.5) circle (.08);
}
\draw[brown] (1.6,3.8) node {$Q'$};
\draw[red] (2.8,3.6) node {$Q$};
\end{tikzpicture}
\caption{A family of paths $F\in\F_\skew$, where $\lambda/\mu=(5,4,4,2,2)/(3,2,1)$ (left), and the tableau $T$ with $\La(T)=F$.}
\label{fig:another_example}
\end{figure}
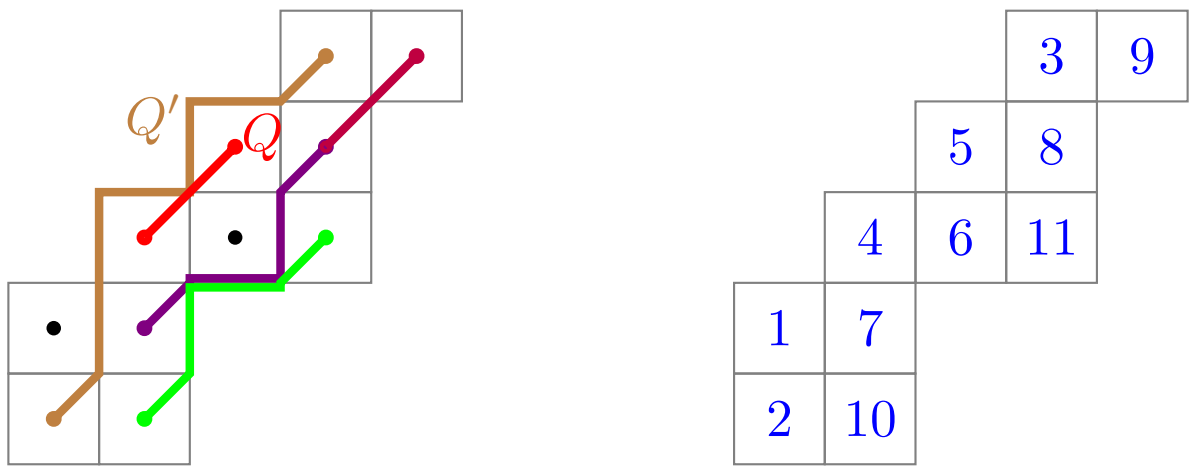

Let $\overline{Q}$ denote the polygonal line formed by the paths in $Q$.
We say that a component $Q$ is left (resp.\ right) of another component $Q'$ if there are points $q\in\overline{Q}\setminus\overline{Q'}$ and $q'\in\overline{Q'}\setminus\overline{Q}$ (not necessarily vertices) such that $q$ and $q'$ are at the same height, and $q$ is strictly left (resp.\ right) of $q'$. Note that, since the paths in $F$ are noncrossing, $Q$ cannot be both left and right of $Q'$. 
In Figure~\ref{fig:F}, $Q_2$ is left of $Q_1$.
Note that it is possible for a component $Q'$ to be left of $Q$ even if $L(Q)=\emptyset$; see Figure~\ref{fig:another_example} for an example of this situation.

Our next goal is to prove that $T$ is uniquely determined by $\La(T)$, and also uniquely determined by $\Ch(T)$. More specifically, we will show that these two statistics provide bijections from $\SYT(\skew)$ to $\F_\skew$. 

We start by describing the inverse of the $\La$ map. Until we prove that it is indeed $\La^{-1}$, we call this map $\Phi$. We label the cells of $\skew$ by rows, as in Definition~\ref{def:labeling_by_rows}, and we say that a cell is {\em smaller} than another if it has a smaller label. In particular, in any subset of cells, the smallest one is the leftmost cell in the top row.

\begin{definition}[The inverse of $\La$] \label{def:La-inv}
    Given a family of paths $F\in\F_\skew$, construct a tableau $\Phi(F)\in\SYT(\skew)$ by filling the cells of $\skew$ with the numbers $1,2,\dots,n$ as follows.
    The components of $F$ will be filled one at a time, in a specific order. The sets $L(Q)$ below consist of unfilled cells only.
    
    Let $Q_1$ be the component containing the smallest cell (with the labeling by rows). If $L(Q_1)\neq \emptyset$, let 
    $Q_2$ be the component containing the smallest cell in $L(Q_1)$. If $L(Q_2)\neq\emptyset$, define $Q_3$ similarly, until eventually we find a component $Q_j$ such that $L(Q_j)=\emptyset$. (This process ends because the smallest cell in $Q_{i+1}$ is in a lower row than the smallest cell in $Q_i$ for all $i$.)
    Fill the cells of $Q_j$ from southwest to northeast. 

    Now that component $Q_j$ has been filled, remove it from $\skew$ and apply the procedure from the previous paragraph recursively to the unfilled cells, until all the cells have been filled. 
\end{definition}

As an example of the construction in Definition~\ref{def:La-inv}, consider the family of paths in Figure~\ref{fig:F}. In the first step of the process, we have $L(Q_1)=\{3,6,7\}$, $L(Q_2)=\{6,11\}$ and $L(Q_3)=\emptyset$, so we start by filling the only cell in $Q_3$ with the entry $1$. The rest of the construction is shown in Figure~\ref{fig:filling}.

\begin{figure}[htb]
\centering
\begin{tikzpicture}[scale=.6]
\drawskewtab{4}{{5,4,1,0}}{{{,},{,,},{1,,,,},{,,,}}}
\drawpath{orange}{(1.5,.5)}{}
\drawpathsmall{red}{(2.5,1.5)}{\e}
\drawpath{brown}{(3.5,1.5)}{\e\n}
\drawpath{violet}{(2.5,.5)}{\e}
\drawpath{green}{(5.5,1.5)}{\n}
\foreach \x in {(1,1),(4,1),(6,3),(7,3)}{
    \fill \x+(-.5,-.5) circle (.08);
}

\begin{scope}[shift={(8,0)}]
\drawskewtab{4}{{5,4,1,0}}{{{,},{,,},{1,,,,},{2,,,}}}
\drawpath{orange}{(1.5,.5)}{}
\drawpathsmall{red}{(2.5,1.5)}{\e}
\drawpath{brown}{(3.5,1.5)}{\e\n}
\drawpath{violet}{(2.5,.5)}{\e}
\drawpath{green}{(5.5,1.5)}{\n}
\foreach \x in {(4,1),(6,3),(7,3)}{
    \fill \x+(-.5,-.5) circle (.08);
}
\end{scope}

\begin{scope}[shift={(16,0)}]
\drawskewtab{4}{{5,4,1,0}}{{{,},{5,,},{1,4,,,},{2,3,,}}}
\drawpath{brown}{(3.5,1.5)}{\e\n}
\drawpath{violet}{(2.5,.5)}{\e}
\drawpath{green}{(5.5,1.5)}{\n}
\foreach \x in {(4,1),(6,3),(7,3)}{
    \fill \x+(-.5,-.5) circle (.08);
}
\end{scope}

\begin{scope}[shift={(0,-5)}]
\drawskewtab{4}{{5,4,1,0}}{{{7,},{5,,},{1,4,6,,},{2,3,,}}}
\drawpath{violet}{(2.5,.5)}{\e}
\drawpath{green}{(5.5,1.5)}{\n}
\foreach \x in {(4,1),(6,3),(7,3)}{
    \fill \x+(-.5,-.5) circle (.08);
}
\end{scope}

\begin{scope}[shift={(8,-5)}]
\drawskewtab{4}{{5,4,1,0}}{{{7,},{5,8,},{1,4,6,,},{2,3,,}}}
\drawpath{violet}{(2.5,.5)}{\e}
\drawpath{green}{(5.5,1.5)}{\n}
\foreach \x in {(4,1),(7,3)}{
    \fill \x+(-.5,-.5) circle (.08);
}
\end{scope}

\begin{scope}[shift={(16,-5)}]
\drawskewtab{4}{{5,4,1,0}}{{{7,},{5,8,},{1,4,6,10,},{2,3,9,}}}
\drawpath{green}{(5.5,1.5)}{\n}
\foreach \x in {(4,1),(7,3)}{
    \fill \x+(-.5,-.5) circle (.08);
}
\end{scope}

\begin{scope}[shift={(0,-10)}]
\drawskewtab{4}{{5,4,1,0}}{{{7,12},{5,8,},{1,4,6,10,11},{2,3,9,}}}
\foreach \x in {(4,1),(7,3)}{
    \fill \x+(-.5,-.5) circle (.08);
}
\end{scope}

\begin{scope}[shift={(8,-10)}]
\drawskewtab{4}{{5,4,1,0}}{{{7,12},{5,8,13},{1,4,6,10,11},{2,3,9,}}}
\foreach \x in {(4,1)}{
    \fill \x+(-.5,-.5) circle (.08);
}
\end{scope}

\begin{scope}[shift={(16,-10)}]
\drawskewtab{4}{{5,4,1,0}}{{{7,12},{5,8,13},{1,4,6,10,11},{2,3,9,14}}}
\end{scope}
\end{tikzpicture}
\caption{Construction of $\La^{-1}(F)=\Phi(F)$ as in Definition~\ref{def:La-inv}, where $F$ is the family of paths from Figure~\ref{fig:F}.}
\label{fig:filling}
\end{figure}

Note that, in the above construction, $L(Q_j)=\emptyset$ and $A(Q_j)=\emptyset$. The latter holds because $A(Q_i)\subseteq L(Q_i)$ at each step of the process. It follows that $\Phi(F)\in\SYT(\skew)$. Next we state another property of this construction.

\begin{lemma}[The {\it left-to-right} property] 
\label{lem:left-to-right}
Let $Q_L$ and $Q_R$ be components of $F\in\F_\skew$ such that $Q_L$ is left of $Q_R$. Then $Q_L$ is filled before $Q_R$ in the construction of $\Phi(F)$ from Definition~\ref{def:La-inv}. 
\end{lemma}

\begin{proof}
    Suppose for contradiction that $Q_R$ is filled before $Q_L$. Consider the step when $Q_R$ is filled. If $L(Q_R)$ contained cells of $Q_L$, then this would contradict the fact that $L(Q_R)=\emptyset$ at the time that $Q_R$ is filled.
    
    Thus, we can assume that $L(Q_R)$ contains no cells of $Q_L$. Since $Q_L$ is left of $Q_R$, the top cell of $Q_L$ must be higher than $Q_R$, and $Q_L$ must be left of the top cell of $Q_R$. Let $Q_1,Q_2,\dots,Q_j=Q_R$ be the components in Definition~\ref{def:La-inv} at the time when $Q_R$ is filled. 
    Since $Q_L$ is left of the top cell of $Q_R=Q_j$, which in turn is left of $Q_{j-1}$ by construction, we deduce that $Q_L$ is left of $Q_{j-1}$. If $Q_L$ is left of the top cell of $Q_{j-1}$, we apply the same argument to deduce that $Q_L$ is also left of $Q_{j-2}$. By repeating this process, we eventually find the largest index $i$ such that $Q_L$ is not left of the top cell of $Q_i$ (such $i$ must exist because $i=1$ satisfies this property by construction).
    At this point, $Q_L$ is left of the top cell of $Q_{i+1}$, but not left of the top cell of $Q_i$. But this contradicts the choice of $Q_{i+1}$ in Definition~\ref{def:La-inv}.
\end{proof}

\begin{theorem}\label{thm:La}
    The map $T\mapsto\La(T)$ is a bijection from $\SYT(\skew)$ to $\F_\skew$. 
\end{theorem}

\begin{proof}
    We already argued that $\La(T)\in\F_\skew$ for any $T\in\SYT(\skew)$. It remains to show that the map $\Phi$ from Definition~\ref{def:La-inv} is the inverse of $\La$. Let $F\in\F_\skew$ and let $T=\Phi(F)$. We will show that $T$ is the only tableau in $\SYT(\skew)$ such that $\La(T)=F$.

    \paragraph{Proof that $\La(T)=F$.} 
    To show that $F\subseteq\La(T)$, take any path $P\in F$, and let $C$ and $C'$ be its endpoints, with $C'$ northeast of $C$. Since the cells of the component containing $P$ are filled from largest to smallest, $C$ is filled immediately before $C'$, so $\{C,C'\}\in\AscC(T)$. Additionally, by Lemma~\ref{lem:left-to-right}, all the cells to the left of $P$ are filled before $C$ and $C'$, and all the cells to the right of $P$ will be filled after. It follows that $P\in\La(T)$.

    To prove that $\La(T)\subseteq F$, let us first show that no ascent of $T$ is created between two different components of $F$. Suppose for contradiction that an ascent is created by filling components $Q$ and $Q'$ right after each other, where the top element of $Q$ is in a lower row than the bottom element of $Q'$. This means that the algorithm chooses to fill $Q'$ right after $Q$ has been filled. But in this case, the algorithm would have chosen $Q'$ instead of $Q$ when both $Q$ and $Q'$ were unfilled, since the presence of $Q$ does not affect the process that determines which component to fill.

    Finally, suppose that $P_i\in\La(T)$. By the previous paragraph, its endpoints $C_i$ and $C_{i+1}$ are in the same component of $F$. Since $P_i$ is determined by the boundary of $T^{<i}$, cells to the left (resp.\ right) of $P_i$ have entries $<i$ (resp.\ $>i+1$), so they were filled before (resp.\ after) $C_i$ and $C_{i+1}$. It follows that $P_i\in F$.

    \paragraph{Proof of uniqueness of $T$.} If $U\in\SYT(\skew)$ satisfies $\La(U)=F$, the cells within each component of $F$ have to be filled consecutively from southwest to northeast. Additionally, we claim that the order in which the components of $F$ are filled must satisfy the left-to-right property from Lemma~\ref{lem:left-to-right}. Indeed, suppose that $Q_L$ and $Q_R$ are components of $F$ such that $Q_L$ is left of $Q_R$. Then either there are cells of $Q_L$ that are left of $Q_R$, or there are cells of $Q_R$ that are right of $Q_L$ (or both). In the first case, the cells of $Q_L$ to the left of $Q_R$ need to be filled before $Q_R$ is filled, in order for $Q_R$ to be a component of $\La(U)=F$. Similarly, in the second case, the cells of $Q_R$ to the right of $Q_L$ need to be unfilled at the time when $Q_L$ is filled, in order for $Q_L$ to be a component of $\La(U)=F$. In both cases, $Q_L$ must be filled before $Q_R$.
    
    Suppose for contradiction that $U\neq T$, where $T$ is the tableau found by our algorithm. Compare the order in which the components of $F$ are filled in $T$ and $U$, and consider the first time where they differ. Let $Q_j$ be the component filled in $T$ at this time, i.e., the one chosen by our algorithm by constructing a sequence $Q_1,Q_2,\dots,Q_j$. Let $Q'$ be the component that is filled immediately before $Q_j$ in $U$. Since $\La(U)=F$, there are no ascents in $U$ between different components of $F$, so the top cell in $Q'$ must be at least as high as the bottom cell in $Q_j$.
    By construction of the sequence $Q_1,Q_2,\dots,Q_j$, any component $Q'$ with this property must have some $Q_i$, where $1\le i\le j$, to its left. But then, the fact that $U$ satisfies the left-to-right property forces $Q_j$ to be filled before $Q'$, a contradiction.
\end{proof}

In the special case that $\skew=\delta_n/\delta_{n-1}$, Theorem~\ref{thm:La} is equivalent to a result of Reading's \cite[Thm.~3.1]{Reading2015}, which describes a bijection from permutations to the noncrossing arc diagrams from Section~\ref{sec:arc_diagrams}. One can think of our inverse map $\La^{-1}$ as a generalization of his map $\delta$. 
One difference to note is that, while Reading proves injectivity of his map by using the fact that a permutation is uniquely determined by its canonical join representation (see Section~\ref{sec:weakorder}  below), we chose to make the proof of Theorem~\ref{thm:La} self-contained without relying on the machinery of semidistributive lattices. 

The analogue of Theorem~\ref{thm:La} for chutes now follows  easily.

\begin{corollary}\label{cor:Ch}
    The map $T\mapsto\Ch(T)$ is a bijection from $\SYT(\skew)$ to $\F_\skew$. 
\end{corollary}

\begin{proof}
    By equation~\eqref{eq:Chutes-Ladders-transpose}, we have $\Ch(T)=\La(T')'$.
    In particular, using Theorem~\ref{thm:La}, the map $T\mapsto\Ch(T)$ is a composition of bijections, hence it is a bijection.
\end{proof}

    It will be useful to have an explicit description of the map $\Ch^{-1}$, which is analogous to the description of $\La^{-1}$ given in Definition~\ref{def:La-inv}, up to conjugation. Specifically, the directions west and east are switched with north and south, respectively.

    Unlike in previous constructions, we now label the cells of $\skew$ by columns, as described in Definition~\ref{def:labeling_by_rows}. In particular, in any subset of cells, the smallest cell (i.e., the one with smallest label) is the top cell in the leftmost column. 

\begin{definition}[The inverse of $\Ch$]\label{def:Ch-inv}
    Given a family of paths $F\in\F_\skew$, construct the tableau $\Ch^{-1}(F)\in\SYT(\skew)$ by filling the cells of $\skew$ with the numbers $1,2,\dots,n$ as follows. The components of $F$ will be filled one at a time, in a specific order. The sets $A(Q)$ below consist of unfilled cells only.
    
    Let $Q_1$ be the component containing the smallest cell (with the labeling by columns). If $A(Q_1)\neq \emptyset$, let 
    $Q_2$ be the component containing the smallest cell in $A(Q_1)$. If $A(Q_2)\neq\emptyset$, define $Q_3$ similarly, until eventually we find a component $Q_j$ such that $A(Q_j)=\emptyset$. Fill the cells of $Q_j$ from northeast to southwest. 

    Now that component $Q_j$ has been filled, remove it from $\skew$ and apply the procedure from the previous paragraph recursively to the unfilled cells, until all the cells have been filled. 
\end{definition}

An example of the construction in Definition~\ref{def:Ch-inv} is given in Figure~\ref{fig:filling-Ch-inv}.

\begin{figure}[htb]
\centering
\begin{tikzpicture}[scale=.6]
\drawskewtab{4}{{5,4,1,0}}{{{,},{,,},{,,,,},{1,,,}}}
\drawpath{orange}{(1.5,.5)}{}
\drawpathsmall{red}{(2.5,1.5)}{\e}
\drawpath{brown}{(3.5,1.5)}{\e\n}
\drawpath{violet}{(2.5,.5)}{\e}
\drawpath{green}{(5.5,1.5)}{\n}
\foreach \x in {(2,2),(4,1),(6,3),(7,3)}{
    \fill \x+(-.5,-.5) circle (.08);
}

\begin{scope}[shift={(8,0)}]
\drawskewtab{4}{{5,4,1,0}}{{{,},{,,},{2,,,,},{1,,,}}}
\drawpath{orange}{(1.5,.5)}{}
\drawpathsmall{red}{(2.5,1.5)}{\e}
\drawpath{brown}{(3.5,1.5)}{\e\n}
\drawpath{violet}{(2.5,.5)}{\e}
\drawpath{green}{(5.5,1.5)}{\n}
\foreach \x in {(4,1),(6,3),(7,3)}{
    \fill \x+(-.5,-.5) circle (.08);
}
\end{scope}

\begin{scope}[shift={(16,0)}]
\drawskewtab{4}{{5,4,1,0}}{{{,},{3,,},{2,4,,,},{1,5,,}}}
\drawpath{brown}{(3.5,1.5)}{\e\n}
\drawpath{violet}{(2.5,.5)}{\e}
\drawpath{green}{(5.5,1.5)}{\n}
\foreach \x in {(4,1),(6,3),(7,3)}{
    \fill \x+(-.5,-.5) circle (.08);
}
\end{scope}

\begin{scope}[shift={(0,-5)}]
\drawskewtab{4}{{5,4,1,0}}{{{6,},{3,,},{2,4,7,,},{1,5,,}}}
\drawpath{violet}{(2.5,.5)}{\e}
\drawpath{green}{(5.5,1.5)}{\n}
\foreach \x in {(4,1),(6,3),(7,3)}{
    \fill \x+(-.5,-.5) circle (.08);
}
\end{scope}

\begin{scope}[shift={(8,-5)}]
\drawskewtab{4}{{5,4,1,0}}{{{6,},{3,,},{2,4,7,8,},{1,5,9,}}}
\drawpath{green}{(5.5,1.5)}{\n}
\foreach \x in {(4,1),(6,3),(7,3)}{
    \fill \x+(-.5,-.5) circle (.08);
}
\end{scope}

\begin{scope}[shift={(16,-5)}]
\drawskewtab{4}{{5,4,1,0}}{{{6,},{3,,},{2,4,7,8,},{1,5,9,10}}}
\drawpath{green}{(5.5,1.5)}{\n}
\foreach \x in {(6,3),(7,3)}{
    \fill \x+(-.5,-.5) circle (.08);
}
\end{scope}

\begin{scope}[shift={(0,-10)}]
\drawskewtab{4}{{5,4,1,0}}{{{6,},{3,11,},{2,4,7,8,},{1,5,9,10}}}
\drawpath{green}{(5.5,1.5)}{\n}
\foreach \x in {(7,3)}{
    \fill \x+(-.5,-.5) circle (.08);
}
\end{scope}

\begin{scope}[shift={(8,-10)}]
\drawskewtab{4}{{5,4,1,0}}{{{6,12},{3,11,},{2,4,7,8,13},{1,5,9,10}}}
\foreach \x in {(7,3)}{
    \fill \x+(-.5,-.5) circle (.08);
}
\end{scope}

\begin{scope}[shift={(16,-10)}]
\drawskewtab{4}{{5,4,1,0}}{{{6,12},{3,11,14},{2,4,7,8,13},{1,5,9,10}}}
\end{scope}
\end{tikzpicture}
\caption{Construction of $\Ch^{-1}(F)$ as in Definition~\ref{def:Ch-inv}, where $F$ is the family of paths from Figure~\ref{fig:F}.}
\label{fig:filling-Ch-inv}
\end{figure}

\subsection{A rowmotion map on standard Young tableaux}
\label{sec:rowmotion}

We are now ready to define rowmotion on standard Young tableaux. Recall the definitions of $\Ch$ and $\La$ given in Definition~\ref{def:chutes&ladders}, and the description of the map $\Ch^{-1}$ given in Definition~\ref{def:Ch-inv}.
Theorem~\ref{thm:La} and Corollary~\ref{cor:Ch} imply that the following map is a bijection.

\begin{definition}[Rowmotion on standard Young tableaux]\label{def:row} Let
$$\row:\SYT(\skew)\to\SYT(\skew)$$
be the bijection defined as $\row=\Ch^{-1}\circ\La$.
\end{definition}

\begin{figure}[htb]
\centering
\begin{tikzpicture}[scale=.6]
\begin{scope}[shift={(0,0)}]
\drawcircpathsmall{orange}{(1.5,.5)}{}
\drawcircpathbig{red}{(2.5,1.5)}{\e}
\drawcircpath{brown}{(3.5,1.5)}{\e\n}
\drawcircpath{violet}{(2.5,.5)}{\e}
\drawcircpath{green}{(5.5,1.5)}{\n}
\drawskewtab{4}{{5,4,1,0}}{{{7,12},{5,8,13},{1,4,6,10,11},{2,3,9,14}}}
\draw[->] (8.5,2) -- node[above]{$\rho$} (9.5,2);
\begin{scope}[shift={(11,0)}]
\drawcircpathsmall{orange}{(1.5,.5)}{}
\drawcircpathbig{red}{(2.5,1.5)}{\e}
\drawcircpath{brown}{(3.5,1.5)}{\e\n}
\drawcircpath{violet}{(2.5,.5)}{\e}
\drawcircpath{green}{(5.5,1.5)}{\n}
    \drawskewtab{4}{{5,4,1,0}}{{{6,12},{3,11,14},{2,4,7,8,13},{1,5,9,10}}}
\end{scope}
\end{scope}

\begin{scope}[shift={(2,-6)}]
\drawcircpath{brown}{(.5,.5)}{\n\n\e\n\e}
\drawcircpath{red}{(1.5,2.5)}{}
\drawcircpathsmall{violet}{(1.5,1.5)}{-- ++(0,.05)\e-- ++(0,-0.05)\n}
\drawcircpathbig{purple}{(3.5,3.5)}{}
\drawcircpath{green}{(1.5,.5)}{\n--++(0,-0.05)\e--++(0,0.05)}
\drawskewtab{5}{{3,2,1,0,0}}{{{3,9},{5,8},{4,6,11},{1,7},{2,10}}}
\draw[->] (6.5,2.5) -- node[above]{$\rho$} (7.5,2.5);
\begin{scope}[shift={(9,0)}]
\drawcircpath{brown}{(.5,.5)}{\n\n\e\n\e}
\drawcircpath{red}{(1.5,2.5)}{}
\drawcircpathsmall{violet}{(1.5,1.5)}{-- ++(0,.05)\e-- ++(0,-0.05)\n}
\drawcircpathbig{purple}{(3.5,3.5)}{}
\drawcircpath{green}{(1.5,.5)}{\n--++(0,-0.05)\e--++(0,0.05)}
    \drawskewtab{5}{{3,2,1,0,0}}{{{2,7},{4,8},{5,6,10},{1,9},{3,11}}}
\end{scope}
\end{scope}

\begin{scope}[shift={(1,-11)}]
\drawcircpath{olive}{(1.5,2.5)}{}
\drawcircpath{red}{(.5,1.5)}{\e\e\n}
\drawcircpath{violet}{(.5,.5)}{\e}
\drawcircpath{orange}{(3.5,2.5)}{}
\drawcircpathsmall{green}{(2.5,.5)}{}
\drawcircpathbig{brown}{(3.5,1.5)}{\n\e}
\drawtab{4}{{{1,2,5,8,10,17},{3,4,6,9,18},{7,11,13,16},{12,14,15,19}}}
\draw[->] (7.5,2) -- node[above]{$\rho$} (8.5,2);
\begin{scope}[shift={(10,0)}]
\drawcircpath{olive}{(1.5,2.5)}{}
\drawcircpath{red}{(.5,1.5)}{\e\e\n}
\drawcircpath{violet}{(.5,.5)}{\e}
\drawcircpath{orange}{(3.5,2.5)}{}
\drawcircpathsmall{green}{(2.5,.5)}{}
\drawcircpathbig{brown}{(3.5,1.5)}{\n\e}
    \drawtab{4}{{{1,3,4,7,13,15},{2,5,6,14,19},{8,9,10,16},{11,12,17,18}}}
\end{scope}
\end{scope}

\end{tikzpicture}

\begin{tikzpicture}[scale=.5]
\drawcircpath{red}{(3.5,3.5)}{\n\e\n\e}
\drawcircpath{violet}{(.5,.5)}{\n\e\n\e\e\n\e\n}
\drawskewtab{8}{{7,6,5,4,3,2,1,0}}{{{1},{3},{6},{4},{2},{7},{8},{5}}}
\draw[->] (8.5,4) -- node[above]{$\rho$} (9.5,4);
\begin{scope}[shift={(11,0)}]
\drawcircpath{red}{(3.5,3.5)}{\n\e\n\e}
\drawcircpath{violet}{(.5,.5)}{\n\e\n\e\e\n\e\n}
\drawskewtab{8}{{7,6,5,4,3,2,1,0}}{{{8},{1},{4},{3},{2},{7},{6},{5}}}
\end{scope}
\end{tikzpicture}
    \caption{Some examples of the rowmotion map $\row$. }
    \label{fig:row_examples}
\end{figure}
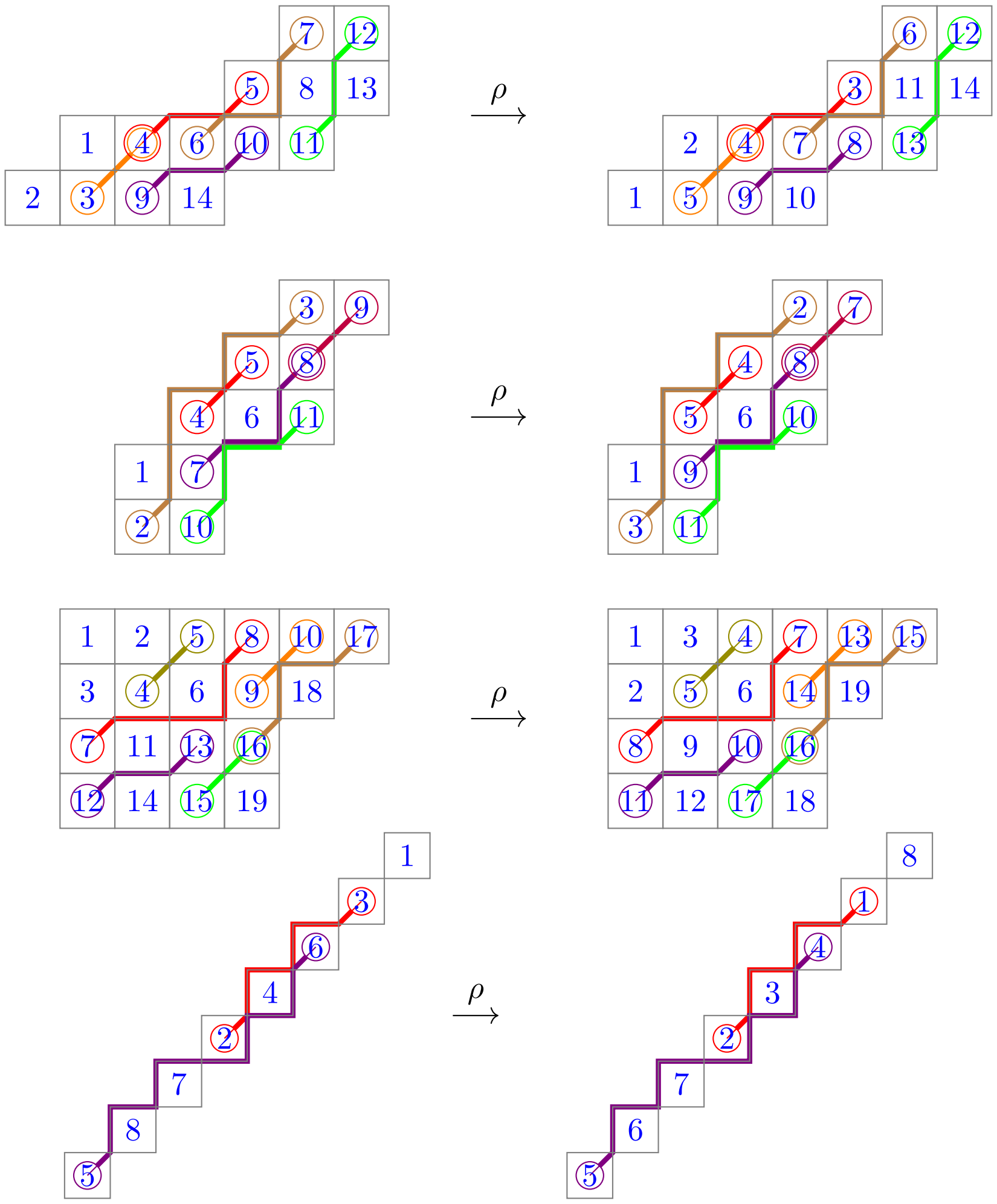

For any $T\in\SYT(\skew)$, its image $\row(T)$ is the tableau whose chutes are precisely the ladders of $T$, that is,
$\La(T)=\Ch(\row(T))$. In particular, $\AscC(T)=\HDesC(\row(T))$ and $\asc(T)=\hdes(\row(T))$.
Figure~\ref{fig:row_examples} shows some examples of the rowmotion map~$\rho$.

\begin{remark}
\label{remark:tworowrho}
When $\skew=(m^2)$, the map $\row$ coincides with order ideal rowmotion on $\A_{m-1}$ (see Section~\ref{sec:two_row}), which is why we use the same notation and terminology for this more general action. See Figure~\ref{fig:2row} for an example. If $T\in\SYT(m^2)$ corresponds to the order ideal $I\in J(\A_{m-1})$, then the ladders of $T$ determine the positions of the minimal elements of $\A_{m-1}\setminus I$,  and the chutes of $T$ determine the positions of the maximal elements of $I$. The fact that these two maps coincide will also follow from Theorem~\ref{thm:kappa}, using Remark~\ref{remark:rowequalsrow}.
\end{remark}

\section{$\kappa$-rowmotion on the semidistributive lattice of standard Young tableaux}
\label{sec:kappa}

In this section, we describe a partial order on $\SYT(\skew)$ and show that the resulting poset is a semidistributive lattice. 
We recall the definition of $\kappa$-rowmotion on semidistributive lattices, and show that, when applied to our lattice, it coincides with the rowmotion map from the previous section. 

We begin with background on semidistributive lattices and the canonical join representation of the weak order on permutations.

\subsection{Semidistributive lattices and canonical join representations}
\label{sec:semidist}
We follow the exposition in \cite[Sec.~3.1]{Barnard2019}; see also \cite[Ch.~3]{EC1} for basic terminology on posets and lattices.
Let $L$ be a finite lattice with meet $\wedge$ and join $\vee$, and partial order denoted by $\leq$. We use $x, y, z$ to denote elements of $L$. We say that $L$ is 
\begin{itemize} 
\item \emph{join-semidistributive} if 
 $x \vee y = x \vee z$ implies $x \vee (y \wedge z) = x \vee y$;
\item \emph{meet-semidistributive} if 
$x \wedge y = x \wedge z$ implies $x \wedge (y \vee z) = x \wedge y$;
\item \emph{semidistributive} if it is both join- and meet-semidistributive.
\end{itemize}

Fix $x\in L$, and let $A\subseteq L$ be such that $x$ can be expressed as the join of all elements in $A$, which we write as $x=\bigvee A:=\bigvee_{a\in A}a$. We say $\bigvee A$ is a \emph{join representation} of $x$. A join representation $\bigvee A$ is \emph{irredundant} if for all strict subsets $A'\subset A$, we have $\bigvee A'<\bigvee A$. Note that a necessary condition for a join representation to be irredundant is that it is an \emph{antichain}, a subset of $L$ with no comparable elements.

Now define a partial order on antichains of $L$. We say that antichain $A$ is \emph{lower} than antichain $B$ if for each $a\in A$ there exists $b\in B$ such that $a\leq b$. Say that $\bigvee A$ 
is a \emph{canonical join representation} of $x$ if it is the unique lowest element of the set of irredundant join representations of $x$. If there is not a unique lowest element, then the canonical join representation does not exist. Each element of a canonical join representation is called a \emph{canonical joinand}. Each canonical joinand is a \emph{join-irreducible} element of $L$, meaning that it cannot be expressed as the join of strictly smaller elements. It is clear that $j$ is join-irreducible if and only if $j$ covers exactly one element of $L$.
Moreover, the canonical join representation of a join-irreducible $j$ is $\{j\}$.

Canonical join representations provide the following alternative characterization of semidistributivity (see \cite[Thm.~16]{Barnard2019}, stating the finite case of~\cite[Thm.~2.24]{FJN1995}): a finite lattice $L$ is join-semidistributive if and only if each element in $L$ has a {canonical join representation}. 

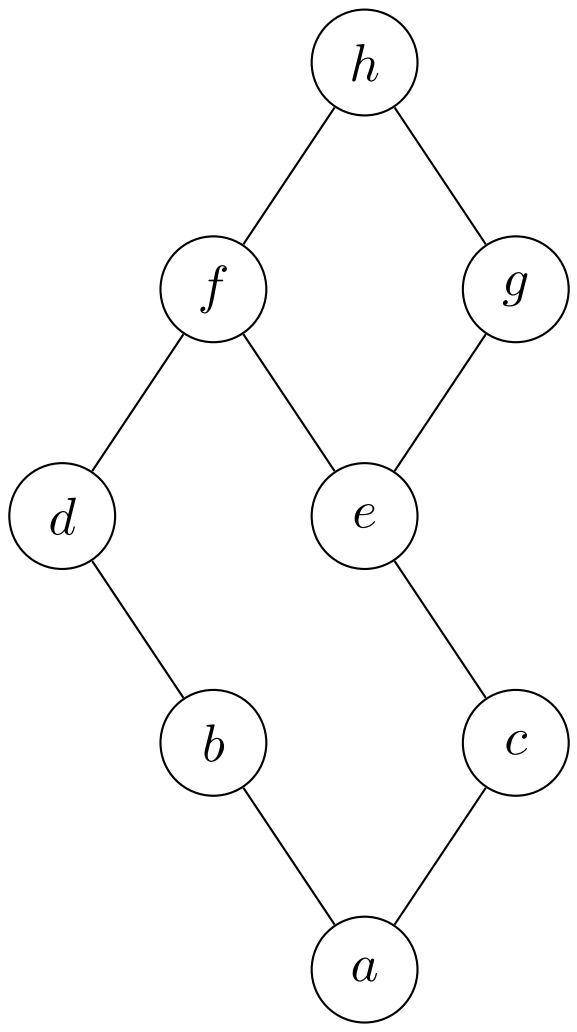
\begin{figure}[htb]
\centering
\begin{tikzpicture}[scale=1]
\tikzset{
    every node/.style={
        circle,
        draw,
        minimum size=7mm,
        inner sep=0pt
    }
}

\node (1) at (1,0) {$a$};
\node (2) at (0,1.5) {$b$};
\node (3) at (2,1.5) {$c$};
\node (4) at (-1,3) {$d$};
\node (5) at (1,3) {$e$};
\node (6) at (0,4.5) {$f$};
\node (7) at (2,4.5) {$g$};
\node (8) at (1,6) {$h$};

\draw (1) -- (2);
\draw (1) -- (3);
\draw (2) -- (4);
\draw (3) -- (5);
\draw (4) -- (6);
\draw (5) -- (6);
\draw (5) -- (7);
\draw (6) -- (8);
\draw (7) -- (8);
\end{tikzpicture}
\caption{A semidistributive lattice.}
\label{fig:poset}
\end{figure}

As an example, consider the semidistributive lattice in Figure~\ref{fig:poset}. 
The elements $b,c,d,e,$ and $g$ are all join-irreducible. The element $f$ can be written as an irredundant join of other elements in several ways: $f=d\vee e=c\vee d=b\vee e=b\vee c$. Of these, the unique lowest antichain is $\{b,c\}$, so the canonical join representation of $f$ is $\bigvee \{b,c\}$. Similarly, the canonical join representation of $h$ is $\bigvee \{b,g\}$. 

\subsection{The canonical join representation of the weak order}
\label{sec:weakorder} 
We now take a slight detour to discuss the weak order on permutations as a semidistributive lattice. 

An {\em inversion} of $\pi\in\S_n$ is a pair of entries $(\pi(i),\pi(j))$ such that $i<j$ and $\pi(i)>\pi(j)$. Let $\Inv(\pi)$ denote the inversion set of $\pi$. The {\em weak order} on $\S_n$ is defined by $\sigma\le\pi$ if $\Inv(\sigma)\subseteq\Inv(\pi)$.
In~\cite{Reading2015}, Reading characterized the canonical join representation of the weak order in terms of permutations and noncrossing arc diagrams. We first review the permutation interpretation. 

\begin{definition}\label{def:oldlambda}
Suppose that $i$ is a descent of $\pi\in \S_n$, and let $b=\pi(i)>\pi(i+1)=a$. Define $\oldlambda(\pi,i)$ to be the permutation whose unique descent is between the values $b>a$
and is such that all values in the interval $(a,b)$ remain on the same side of the descent as they are in $\pi$. More formally, let $(a,b)_L$ (resp.\ $(a,b)_R$) denote the values in the interval $(a,b)$ that appear to the left (resp.\ right) of $a$ and $b$ in $\pi$, written in increasing order. Then $\oldlambda(\pi,i)$ is the permutation $1,\ldots,a-1,(a,b)_L,b,a,(a,b)_R,b+1,\ldots,n$.
\end{definition}

\begin{theorem}[\protect{\cite[Thm.~2.4]{Reading2015}}]
\label{thm:Reading_cjr_perm}
     The canonical join representation of $\pi\in\S_n$ in the weak order  is $$\pi=\bigvee \{\oldlambda(\pi,i): i\in\Des(\pi)\}.$$
\end{theorem}

This canonical join representation can be more easily visualized by interpreting permutations as noncrossing arc diagrams, as described in Section~\ref{sec:arc_diagrams}. 
Recall that for each descent $b>a$ of $\pi$, the arc diagram of $\pi$ has an arc between dots $a$ and $b$ that passes to the right of the dots with labels in $(a,b)_L$ and to the left of the dots with labels in $(a,b)_R$.
For each arc in the arc diagram of $\pi$, one can consider a new arc diagram having only this one arc. The arc diagrams obtained in this way correspond to the permutations $\oldlambda(\pi,i)$ in the canonical join representation of~$\pi$.

\subsection{A partial order on standard Young tableaux}\label{sec:SYTposet}
In this section we describe a partial order on the set $\SYT(\skew)$, which already appeared in work of Bj\"orner and Wachs \cite{Bjorner_Wachs}.

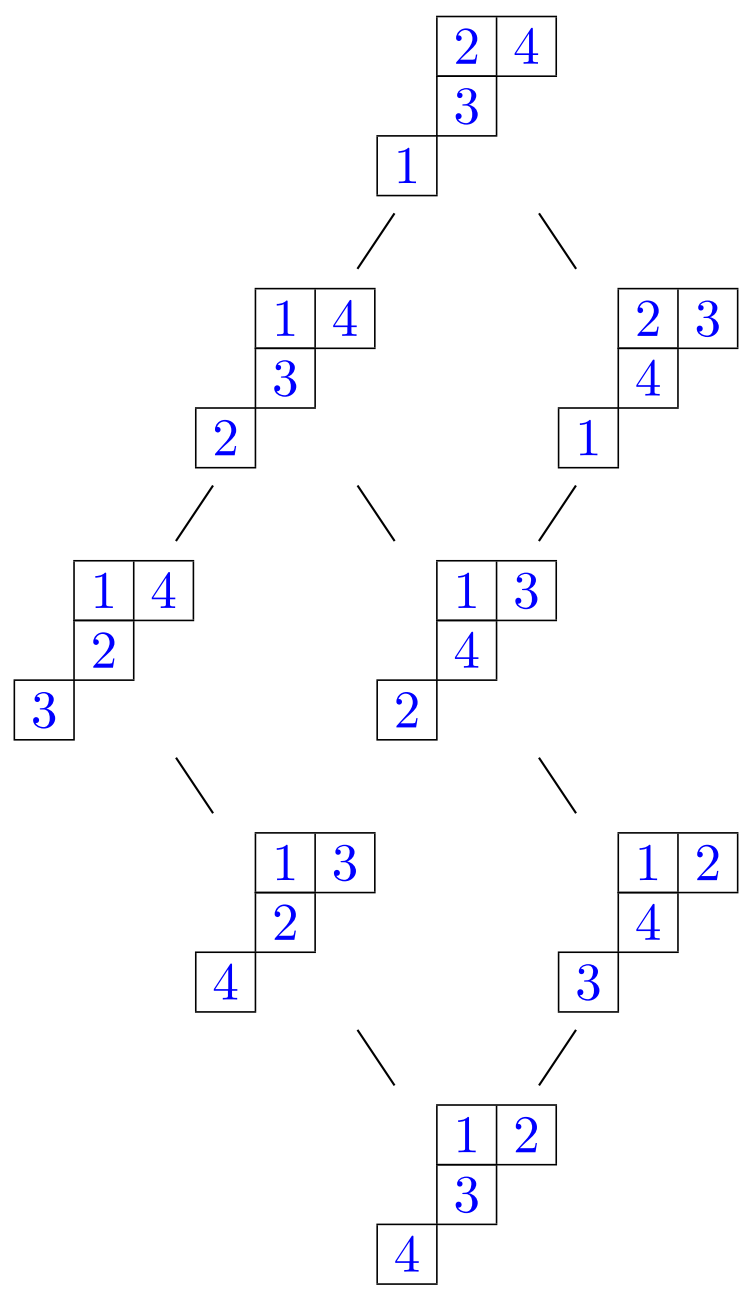
\begin{figure}[htb]
\centering
\begin{tikzpicture}[scale=1.2]
\node (1) at (1,0) {$\young(:\one\two,:\three,\four)$};
\node (2) at (0,1.5) {$\young(:\one\three,:\two,\four)$};
\node (3) at (2,1.5) {$\young(:\one\two,:\four,\three)$};
\node (4) at (-1,3) {$\young(:\one\four,:\two,\three)$};
\node (5) at (1,3) {$\young(:\one\three,:\four,\two)$};
\node (6) at (0,4.5) 
{$\young(:\one\four,:\three,\two)$};
\node (7) at (2,4.5) 
{$\young(:\two\three,:\four,\one)$};
\node (8) at (1,6) {$\young(:\two\four,:\three,\one)$};

\draw (1) -- node[left] {} (2);
\draw (1) -- node[right] {} (3);
\draw (2) -- node[left] {} (4);
\draw (3) -- node[right] {} (5);
\draw (4) -- node[left] {} (6);
\draw (5) -- node[left] {} (6);
\draw (5) -- node[right] {} (7);
\draw (6) -- node[left] {} (8);
\draw (7) -- node[right] {} (8);

\end{tikzpicture}
\caption{The semidistributive lattice $\T_{(3,2,1) / (1,1)}$. }
\label{fig:T321_11}
\end{figure}

We denote by $\T_\skew$ the poset whose elements are tableaux in $\SYT(\skew)$, and whose 
covering relations are given by $T\lessdot U$ if $U$ is obtained from $T$ by swapping entries $i$ and $i+1$ for some $1\le i\le n-1$, where $i$ is northeast of $i+1$ (i.e.~a nontrivial Bender--Knuth involution). In other words, $i$ is a high descent of $T$, whereas $i$ is an ascent of $U$. See Figure~\ref{fig:T321_11} for an example.

The minimal element of $\T_\skew$ is the tableau $T_R$ defined in Section~\ref{sec:SYT}, since it is the only tableau with no ascents. 
The maximal element is $T_C$, which is the only tableau with no high descents.  

\begin{definition}\label{def:pi_T}
To each $T\in\SYT(\skew)$, we associate a permutation $\pi_T\in\S_n$ as follows.
First, label the cells of $\skew$ by rows, as in Definition~\ref{def:labeling_by_rows}. Let $\pi_T$ be the permutation obtained by listing the labels of the cells containing the entries $1,2,\dots,n$ of $T$, in this order. Specifically, $\pi_T(i)=T_R(T^{-1}(i))$ for all $1\le i\le n$. 
\end{definition}

Table~\ref{tab:pi_T} shows two examples.
Note that $\pi^{-1}_T$ is the {\em reading word} of $T$, obtained by reading its entries row by row from left to right in the usual reading order.

\renewcommand{\arraystretch}{1.3}
\begin{table}[htb]
\[
\begin{array}{c||c|c}
T & \young(\one\three\four\eight,\two\six\nine,\five\seven)  & \young(:\two\four\six,:\three\seven,\one\five) \\[15pt] \hline
\pi_T & 152386947 &  6142735 \\ \hline
\pi_T^{-1} & \textcolor{blue}{134826957}  & \textcolor{blue}{2463715}
\end{array}
\]
    \caption{Two examples of the map $T\mapsto\pi_T$.}
    \label{tab:pi_T}
\end{table}

The permutation $\pi_{T_R}$ is the identity permutation $\id:=12\dots n$. The permutation $\pi_{T_C}$ will be denoted by $\sigma_{\lambda/\mu}$. For example, if $\lambda=(4,3,2)$, then $\sigma_{\lambda}=158269374$, and if $\lambda/\mu=(4,3,2)/(1,1,0)$, then
$\sigma_{\lambda/\mu}=6147253$.

The following result follows from \cite[Thm.~7.2]{Bjorner_Wachs}. We include a proof here for the sake of self-containment.

\begin{proposition}
\label{prop:interval}
    For any skew shape $\skew$ with $n$ cells, the map $T\to\pi_T$ is a poset isomorphism between $\T_\skew$ and the interval $[\id,\sigma_\skew]$ in the weak order on $\S_n$. In addition,
    \begin{equation}\label{eq:AscDes} \Asc(T)=\Des(\pi_T).\end{equation}
\end{proposition}

\begin{proof}
The map $T\to\pi_T$ is clearly injective, since $T$ is uniquely determined by $\pi_T$. Thus, it is a bijection between $\SYT(\skew)$ and the set $\mathcal{I}_\skew:=\{\pi_T:T\in\SYT(\skew)\}$. Let us characterize the permutations in this set. Recall that $\pi_T$ is obtained by reading the labels (in the labeling by rows, which is given by the entries of $T_R$) of the cells containing the entries $1,2,\dots,n$ of $T$.

The fact that the entries of $T$ increase along rows (from left to right) and along columns (from top to bottom) translates into the requirement that if $a$ and $b$ are in the same row or column of $T_R$ and $a<b$, then $a$ appears before $b$ in $\pi_T$. By transitivity, this is equivalent to requiring that if $b$ is weakly southeast of $a$ in $T_R$, then $a$ appears before $b$ in $\pi_T$.
Thus, $\mathcal{I}_\skew$ is the set of permutations $\pi\in\S_n$ such that $(b,a)\notin\Inv(\pi)$ for any pair of entries such that $b$ is weakly southeast of $a$ in $T_R$. Equivalently, for every $(b,a)\in\Inv(\pi)$, entry $a$ must be (strictly) northeast of $b$ in $T_R$. By construction, any such pair $(b,a)$ is an inversion of the permutation $\sigma_\skew$. Thus, 
$$\mathcal{I}_\skew=\{\pi\in\S_n: \Inv(\pi)\subseteq\Inv(\sigma_\skew)\}=\{\pi\in\S_n: \pi\le\sigma_\skew \text{ in the weak order}\}=[\id,\sigma_\skew].$$

This bijection is a poset isomorphism because the covering relations are the same in $\T_\skew$ as in the interval $[\id,\sigma_\skew]$ in the weak order. Indeed, $T\lessdot U$ if $U$ is obtained from $T$ by swapping entries $i$ and $i+1$, where $i$ is strictly northeast of $i+1$. In the corresponding permutations, $\pi_{U}$ is obtained from $\pi_{T}$ by swapping entries in positions $i$ and $i+1$, where $\pi_T(i)<\pi_T(i+1)$.

Finally, recall that $i\in\Asc(T)$ if $i$ is in a lower row than $i+1$ in $T$. This happens if and only if the cell containing $i$ has a bigger label (in the labeling by rows) than the cell containing $i+1$, that is, $\pi_T(i)>\pi_T(i+1)$. 
\end{proof}

Recall the notation $\delta_n=(n,n-1,\dots,1)$ for the staircase shape.

\begin{corollary}
\label{cor:disjoint}
    The poset $\T_{\delta_n/\delta_{n-1}}$ is isomorphic to the weak order on $\S_n$. 
\end{corollary}
\begin{proof}
    The skew shape $\delta_n/\delta_{n-1}$ consists of $n$ cells, arranged in a southwest-northeast diagonal line. Thus, $\sigma_{\delta_n/\delta_{n-1}}$ is the permutation $n(n-1)\ldots 21$, which is the top element in the weak order on $\S_n$. By Proposition~\ref{prop:interval}, $\T_{\delta_n/\delta_{n-1}}$ is isomorphic to the full weak order.
\end{proof}

An interval of the form $[\id,\pi]$ in the weak order will be called a lower interval. By Proposition~\ref{prop:interval}, $\T_{\skew}$ is always isomorphic to a lower interval in the weak order. It is natural to ask whether every lower interval of the weak order arises in this way; equivalently, whether every permutation in $\S_n$ arises as $\sigma_{\skew}$ for some skew shape $\skew$. It is easy to see that this is not the case in general. For instance, the permutation $1243$ is not of the form $\sigma_{\skew}$ for any $\skew$.
It would be interesting to characterize the permutations that arise in this way; equivalently, which lower intervals of the weak order are isomorphic to $\T_\skew$ for some skew shape.

The following is another consequence of Proposition~\ref{prop:interval}.

\begin{corollary}
\label{cor:semidistr}
    For any skew shape $\skew$, the poset $\T_\skew$ is a semidistributive lattice.
\end{corollary}
\begin{proof}
    Proposition~\ref{prop:interval} shows that $\T_\skew$ is isomorphic to a lower interval in the weak order on $\S_n$, which is a semidistributive lattice by \cite[Thm.~3]{DC94} (see also \cite[Sec.~2]{Reading2015}). Any interval in a semidistributive lattice is itself semidistributive, by the definition of semidistributivity in terms of meets and joins (as in Section~\ref{sec:semidist}). 
    Thus, $\T_\skew$ is a semidistributive lattice.  
\end{proof}

As an example, the poset $\T_{(3,2,1)/(1,1)}$ from Figure~\ref{fig:T321_11} is a  semidistributive lattice, and in fact is isomorphic to the one in Figure~\ref{fig:poset}.

Our next goal is to characterize the skew shapes $\skew$ for which $\T_\skew$ has the stronger property of being distributive. A lattice $L$ is \emph{distributive} if, for all
$x, y, z\in L$, we have
 $$x \vee (y \wedge z) = (x \vee y) \wedge (x \vee z) \quad \text{and} \quad x \wedge (y \vee z) = (x \wedge y) \vee (x \wedge z).$$
We will need the following lemmas.

 \begin{lemma}
 \label{lem:dual}
    The posets $\T_{\skew}$ and $\T_{\lambda'/\mu'}$ are dual to each other.
\end{lemma}

\begin{proof}
    Suppose that $T\lessdot U$ in $\T_{\skew}$, where $i$ is northeast of $i+1$ in $T$, and $U$ is obtained from $T$ by swapping these two entries. Then $U'\lessdot T'$ in $\T_{\lambda'/\mu'}$, because $i$ is northeast of $i+1$ in $U'$, and $T'$ is obtained from $U'$ by swapping the entries $i$ and $i+1$.
\end{proof}

\begin{lemma}\label{lem:3rows-cols}
    Let $\skew$ be a skew shape. Then $\skew$ contains three cells $\{C_1,C_2,C_3\}$ such that $C_{i}$ is strictly northeast of $C_{i+1}$ for $i\in\{1,2\}$ if and only if $\skew$ has at least three rows, it has at least three columns, and it is not contained in a hook or a reverse hook.
\end{lemma}

\begin{proof}
    If $\skew$ contains three cells $\{C_1,C_2,C_3\}$ satisfying the above conditions, then these cells are in three different rows, in three different columns, and are not contained in a hook or reverse hook.

    To prove the converse, we will repeatedly use the following property of skew shapes, which we call the rectangle property: if a skew shape contains two cells $(i_1,j_1)$ and $(i_2,j_2)$ with $i_1<i_2$ and $j_1<j_2$, then it must also contain the cells $(i_1,j_2)$ and $(i_2,j_1)$.

    Now let $\skew$ be a skew shape with at least three rows, at least three columns, and not contained in a hook or a reverse hook. 
    Let $i_1$ and $i_3$ be the top and the bottom rows of $\skew$, and let $j_1$ and $j_3$ be the leftmost and rightmost columns of $\skew$.
    
    Since $\skew$ has cells in row $i_1$ and in column $j_3$, the rectangle property implies that $\skew$ contains the cell $C_1=(i_1,j_3)$. Similarly, since $\skew$ has cells in row $i_3$ and in column $j_1$, it must contain the cell $C_3=(i_3,j_1)$. If $\skew$ also contains a cell $C_2=(i,j)$ with $i_1<i<i_3$ and $j_1<j<j_3$, then we have found a triple $\{C_1,C_2,C_3\}$ as desired. 
    
    In the rest of the proof, we will assume that $\skew$ contains no such cell. Since $\skew$ has at least 3 rows, it must have a cell in some row $i_2$ with $i_1<i_2<i_3$. By the above assumption, this cell must be $(i_2,j_1)$ or $(i_2,j_3)$.
    Similarly, since $\skew$ has at least 3 columns, there is some column $j_2$ with $j_1<j_2<j_3$ so that $\skew$ contains $(i_1,j_2)$ or $(i_3,j_2)$. Let us look at the four possibilities depending on which of these cells are in $\skew$.

    If $\skew$ contains $(i_2,j_1)$ and $(i_3,j_2)$, then it would contain $(i_2,j_2)$ by the rectangle property, which is a contradiction. If $\skew$ contains $(i_2,j_3)$ and $(i_1,j_2)$, we reach the same conclusion.

    If $\skew$ contains $(i_2,j_1)$ and $(i_1,j_2)$, then the fact that $\skew$ is not contained in a hook forces it to have some cell $(i,j)$ with $i>i_1$ and $j>j_1$. By the assumption, we must have $i=i_3$ or $j=j_3$. By symmetry, suppose that $j=j_3$. Then the cells $(i_1,j_2)$ and $(i,j_3)$, together with the rectangle property, force $\skew$ to contain the cell $(i,j_2)$, which can only happen if $i=i_3$. But now the cells $(i_3,j_2)$ and $(i_2,j_1)$ force $\skew$ to also contain the cell $(i_2,j_2)$, which is a contradiction. 

    Finally, if $\skew$ contains $(i_3,j_2)$ and $(i_2,j_3)$, then the fact that $\skew$ is not contained in a reverse hook forces it to have a cell $(i,j)$ with $i<i_3$ and $j<j_3$, and we can use a similar argument to reach a contradiction.   
\end{proof}

\begin{theorem}
\label{thm:dist}
    The lattice $\T_{\skew}$ is distributive if and only if $\skew$ has at most two rows or at most two columns, or is contained in a hook or a reverse hook.
\end{theorem}
\begin{proof}
    Suppose $\skew$ has at least three rows and  at least three columns and is not contained in a hook or reverse hook. By Lemma~\ref{lem:3rows-cols}, $\skew$ contains a set of three cells $\mathcal{C}=\{C_1,C_2,C_3\}$ where $C_{i}$ is strictly northeast of $C_{i+1}$ for $i\in\{1,2\}$. 
    Given $\tau\in\S_3$, let $T^{\mathcal{C}}_{\tau}$ denote the tableau in $\SYT(\skew)$ constructed as follows. Consider the skew shape consisting of cells weakly northwest of any cell of $\mathcal{C}$, excluding the cells of $\mathcal{C}$, and suppose that this shape has $k$ cells. 
    Fill these cells with $1,2,\ldots,k$ in reading order, so that all the entries in each row are smaller than all the entries in the next row. Next, for each $1\leq i\leq 3$, place $k+\tau(i)$ in cell $C_i$. Finally, place $k+4,\ldots,|\skew|$ in the remaining cells again in reading order.  
    
    Consider the interval $[T^{\mathcal{C}}_{123},T^{\mathcal{C}}_{321}]$. The standard Young tableaux in this interval are obtained by permuting the entries in $\mathcal{C}$. Thus, this interval is isomorphic to $\T_{(3,2,1)/(2,1)}$. By Corollary~\ref{cor:disjoint}, this poset is isomorphic to the weak order on $\S_3$, which is not distributive. 
    Hence, $\T_\skew$ has an interval that is not distributive, which implies that it is itself not distributive. 

    Now suppose $\skew$ has at most two rows or at most two columns, or is contained in a hook or a reverse hook. We start with the case where $\skew$ has at most two rows, so $\skew\subseteq (m,m)$ for some $m$. The poset $\T_{(m,m)}$ is a distributive lattice, isomorphic to the lattice of order ideals of $\A_{m-1}$ ordered by containment~\cite[Thm.~4.12]{prorow}.
    In particular, $\T_{\skew}$ is isomorphic to an interval of $\T_{(m,m)}$, so it is also a distributive lattice.  
    
    If $\skew$ has at most two columns, its transpose ${\lambda'/\mu'}$ has at most two rows. By Lemma~\ref{lem:dual}, the posets $\T_{\skew}$ and $\T_{\lambda'/\mu'}$ are dual to each other, so the distributivity of $\T_{\skew}$ follows from the distributivity of $\T_{\lambda'/\mu'}$.
        
    If $\skew$ is contained in a hook, first observe that standard Young tableaux of straight shape $\nu$ always have a $1$ in the upper left cell, and so $\T_{\nu}$ is isomorphic to $\T_{\nu/(1)}$. In particular, if $\nu$ is a hook, say $\nu=(a+1,1^b)$, then $\T_{\nu}\cong\T_{\nu/(1)}\cong \T_{(a+b,b) /(b)}$, which is a distributive lattice since $(a+b,b) /(b)$ has at most two rows.  
    For any shape $\skew$ contained in a hook, we similarly deduce that $\T_{\skew}\cong \T_{(a+b,b) /(b)}$ for some $a$ and $b$.
    
    An analogous argument applies to shapes contained in a reverse hook, using that standard Young tableaux of reverse hook shape always have their largest entry in the bottom right corner.
\end{proof}

An alternative approach to the above proof would be to use a result of Stembridge \cite[Thm.~2.2]{Stembridge96}, which states that the interval $[\id, \pi]$ in the weak order is a distributive lattice if and only if $\pi$ is \emph{fully commutative}, 
which is equivalent to $\pi$ avoiding the pattern $321$. It is not hard to see that a skew shape $\skew$ contains three cells $\{C_1,C_2,C_3\}$ as in Lemma~\ref{lem:3rows-cols} if and only if $\sigma_{\skew}$ contains the pattern $321$. 

An interesting problem is to determine for which skew shapes the lattice $\T_\skew$ is rank-symmetric. Rank-symmetry of weak order intervals has been studied; for instance, \cite{Wei12} shows that if $\pi$ is separable (i.e., avoids $3142$ and $2413$), then $[\id,\pi]$ is rank-symmetric and rank unimodal. No general characterization is known, but it was shown in \cite{Denoncourt} that at least $\lfloor \frac{n}{2}\rfloor !$ of the 
intervals of the form $[\id,\pi]$ in the weak order on $\S_n$ are rank-symmetric.  The smallest shape $\skew$ for which $\T_{\skew}$ is not rank-symmetric is $(3,2) / (1)$ (where $\sigma_{(3,2) / (1)}=3142$), and the smallest such straight shape is $(3,2)$ (where $\sigma_{(3,2)}=14253$).

\subsection{$\kappa$-rowmotion on semidistributive lattices}
\label{sec:kapparow}
A rowmotion operation $\kappa:L\rightarrow L$ on an arbitrary finite semidistributive lattice $L$ was defined in \cite{Barnard2019}. After discussing this construction, the main result of this section will be Theorem~\ref{thm:kappa}, showing that when this abstract definition is applied to $\SYT(\skew)$, it matches our definition of $\rho$ in terms of chutes and ladders from Section~\ref{sec:rowmotion-chutes&ladders}. 

For any cover relation $y\lessdot x$, we will denote the corresponding edge of the Hasse diagram by $yx$.
Let $L$ be a finite semidistributive lattice, and let $x\in L$. As shown in \cite[Lem.~19, part 1]{Barnard2019}, there is a bijection between the set of elements covered by $x$ and the set of canonical joinands of $x$, which sends an element $y\lessdot x$ to the unique canonical joinand $j$ such that $y \vee j = x$. 
This bijection is used to label the edges of the Hasse diagram of $L$ as follows.

\begin{definition}
\label{def:labels}
For each covering relation $y\lessdot x$ in $L$, the label of the edge $yx$, denoted by $\ell(yx)$, is the unique canonical joinand $j$ such that $y \vee j = x$. 
\end{definition}

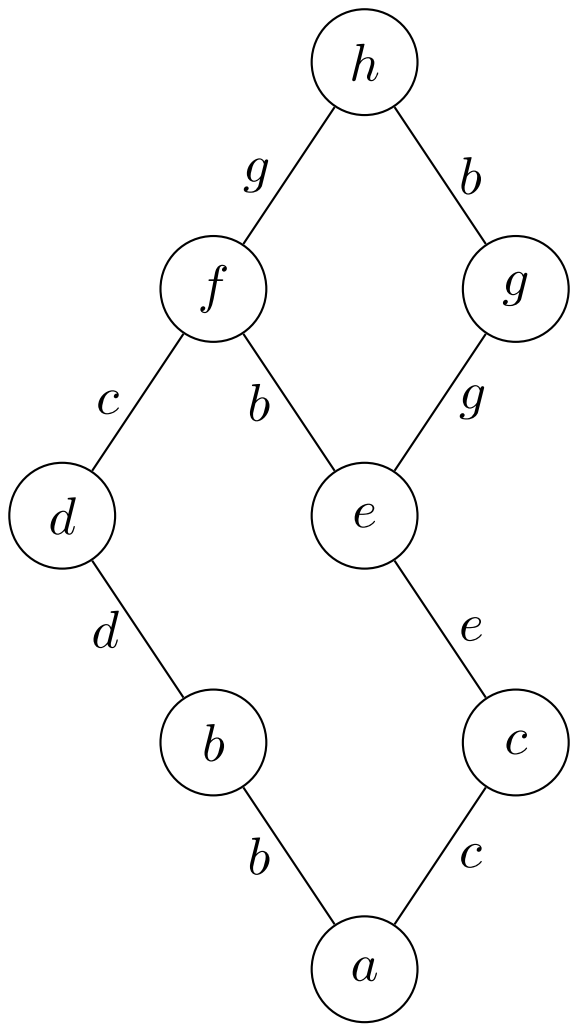
\begin{figure}[htb]
\centering
\begin{tikzpicture}[scale=1]
\tikzset{
    vertex/.style={
        circle,
        draw,
        minimum size=7mm,
        inner sep=0pt
    }
}

\node[vertex] (1) at (1,0) {$a$};
\node[vertex] (2) at (0,1.5) {$b$};
\node[vertex] (3) at (2,1.5) {$c$};
\node[vertex] (4) at (-1,3) {$d$};
\node[vertex] (5) at (1,3) {$e$};
\node[vertex] (6) at (0,4.5) {$f$};
\node[vertex] (7) at (2,4.5) {$g$};
\node[vertex] (8) at (1,6) {$h$};

\draw (1) -- node[left] {$b$} (2);
\draw (1) -- node[right] {$c$} (3);
\draw (2) -- node[left] {$d$} (4);
\draw (3) -- node[right] {$e$} (5);
\draw (4) -- node[left] {$c$} (6);
\draw (5) -- node[left] {$b$} (6);
\draw (5) -- node[right] {$g$} (7);
\draw (6) -- node[left] {$g$} (8);
\draw (7) -- node[right] {$b$} (8);

\end{tikzpicture}
\caption{An example of the edge labels given by Definition~\ref{def:labels}.}
\label{fig:poset_labels}
\end{figure}

Figure~\ref{fig:poset_labels} shows the edge labels for the semidistributive lattice from Figure~\ref{fig:poset}. 

For $x\in L$, define the set of \emph{lower cover labels} of $x$ as $\lowlabels(x)=\{\ell(yx) \ | \ y\lessdot x\}$ and the set of \emph{upper cover labels} of $x$ as $\uplabels(x)=\{\ell(xy) \ | \ x\lessdot y\}$. For example, in the lattice from Figure~\ref{fig:poset_labels}, we have $\lowlabels(f)=\{b,c\}$ and $\uplabels(f)=\{g\}$.

The following property of semidistributive lattices is due to Barnard.

\begin{theorem}[\protect{\cite[Cor.~5]{Barnard2019}}]
\label{thm:Barnard}
Let $L$ be a finite semidistributive lattice and $x\in L$. There is a unique element $y\in L$ with $\uplabels(y)=\lowlabels(x)$.
\end{theorem}

The following is equivalent to the definition of $\kappa$-rowmotion given in~\cite{Barnard2019}.

\begin{definition}[$\kappa$-rowmotion on a semidistributive lattice]
\label{def:kappa}
Let $L$ be a finite semidistributive lattice. The \emph{$\kappa$-rowmotion} map $\kappa:L\rightarrow L$ is defined by letting $\kappa(x)=y$, where $y$ is the unique element in $L$ such that $\uplabels(y)=\lowlabels(x)$. 
\end{definition}

To compute $\kappa(x)$, we look at the set of lower cover labels of $x$ and find the unique element in $L$ that has this set as its set of upper cover labels.
As noted by Barnard~\cite{Barnard2019}, it follows from~\cite{FJN1995} that the map $\kappa:L\rightarrow L$ is a bijection.

As an example, consider the semidistributive lattice in Figure~\ref{fig:poset_labels}. To compute $\kappa(f)$, we note that $\lowlabels(f)=\{b,c\}=\uplabels(a)$. The action of $\kappa$ on all the elements of this lattice is as follows:
$$f\xrightarrow{\kappa}a\xrightarrow{\kappa}h\xrightarrow{\kappa}e\xrightarrow{\kappa}c\xrightarrow{\kappa}d\xrightarrow{\kappa}b\xrightarrow{\kappa}g\xrightarrow{\kappa}f.$$
In this case, all the elements of the lattice are in the same $\kappa$-orbit, but in general there may be more than one orbit. 

\begin{remark}
\label{remark:rowequalsrow}
    If $L$ is a distributive lattice, then it is isomorphic to the lattice of order ideals $J(P)$ for some poset $P$ (its subposet of join-irreducible elements), by the fundamental theorem of finite distributive lattices. In this case, $\kappa$-rowmotion on $L$ is equivalent to order ideal rowmotion on $P$; see \cite[Rmk.\ 32 and Prop.\ 33]{Barnard2019}.
\end{remark}

\subsection{$\kappa$-rowmotion on standard Young tableaux}
We now state our main theorem relating the lattice-theoretic $\kappa$-rowmotion map from Definition~\ref{def:kappa} on the semidistributive lattice $\T_{\skew}$ to the rowmotion map $\rho$ from Definition~\ref{def:row}, justifying our use of the same word for both maps. The rest of the subsection is devoted to its proof.

\begin{theorem}\label{thm:kappa}
    For any skew shape $\skew$ and any $T\in\SYT(\skew)$, we have $$\rho(T)=\kappa(T).$$
\end{theorem}

Recall Reading's characterization of the canonical join representation of the weak order in terms of permutations~\cite{Reading2015}, stated as Theorem~\ref{thm:Reading_cjr_perm} above. We use this to characterize the canonical join representation of elements of $\T_\skew$ in terms of chutes and ladders by relating the permutations $\oldlambda(\pi,i)$ from Definition~\ref{def:oldlambda}
to tableaux. 

\begin{lemma}
 \label{lem:oldlambdatab}
     For any tableau $T\in\SYT(\skew)$ and $i\in\Asc(T)$, we have  
     \[\oldlambda(\pi_T,i)=\pi_{\La^{-1}(\{P_i\})},\] 
     where $P_i$ is the ladder in $\La(T)$ from the cell containing $i$ to the cell containing $i+1$. 
\end{lemma}
 \begin{proof}
 Suppose that, in the labeling by rows of $\skew$, the cells containing $i+1$ and $i$ in $T$ have labels $a$ and $b$, respectively. We refer to these cells as cells $a$ and $b$. 
 
     By Theorem~\ref{thm:La}, there is a unique $U\in\SYT(\skew)$ such that $\La(U)=\{P_i\}$.  The tableau $U$ can be obtained by filling the cells of $\skew$ as follows. Start by filling the cells in reading order (i.e., row by row from left to right and top to bottom) with the numbers $1,\ldots,a-1$, until reaching the cell just before cell $a$. Let $k$ be the number of unfilled cells to the left of the path $P_i$, and fill them with the numbers $a,a+1,\ldots,a+k-1$ again in reading order. Fill cell $b$ with $a+k$ and cell $a$ with $a+k+1$. Finally, fill the remaining cells with the numbers $a+k+2,a+k+3,\ldots,n$ in reading order. 
     By construction, the permutation $\pi_{U}=\pi_{\La^{-1}(\{P_i\})}$ is precisely $\oldlambda(\pi_T,i)$. 
 \end{proof}

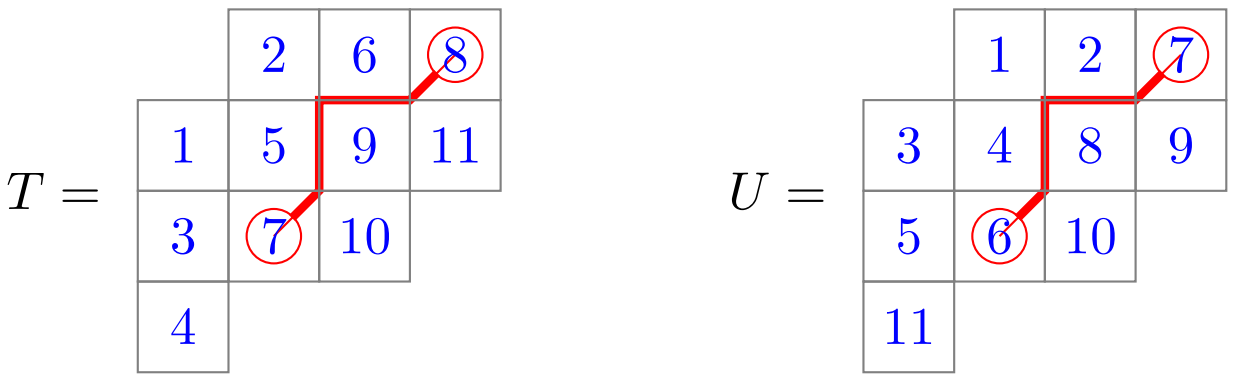
\begin{figure}[htb]
\centering
\begin{tikzpicture}[scale=.6]
\draw (-.2,2) node[left] {$T=$};
\drawcircpath{red}{(1.5,1.5)}{\n\e}
\drawskewtab{4}{{1,0,0,0}}{{{2,6,8},{1,5,9,11},{3,7,10},{4}}}
\begin{scope}[shift={(8,0)}]
    \draw (-.2,2) node[left] {$U=$};
\drawcircpath{red}{(1.5,1.5)}{\n\e}
\drawskewtab{4}{{1,0,0,0}}{{{1,2,7},{3,4,8,9},{5,6,10},{11}}}
\end{scope}

\end{tikzpicture}
    \caption{An example of the construction in Lemma~\ref{lem:oldlambdatab}.}
    \label{fig:oldlambdatab}
\end{figure}

For example, let $T$ be the tableau on the left of Figure~\ref{fig:oldlambdatab}, which has $\pi_T=4\,1\,8\,11\,5\,2\,9\,3\,6\,10\,7$, and let $i=7$, so that $a=3$ and $b=9$. The path $P_i$ is drawn in red. Then $U=\La^{-1}(\{P_i\})$ is the tableau on the right, which has $\pi_U=1\,2\,4\,5\,8\,9\,3\,6\,7\,10\,11=\oldlambda(\pi_T,7)$.

 \begin{lemma}
 \label{lem:cjr_tab}
 The canonical join representation of $T\in\SYT(\skew)$ in the lattice $\T_\skew$ is $$T=\bigvee \{\La^{-1}(\{P\}): P\in\La(T)\}.$$
\end{lemma}
\begin{proof}
By Proposition~\ref{prop:interval}, the map $T\mapsto\pi_T$ is a poset isomorphism between $\T_\skew$ and a lower interval in the weak order of $\S_n$. The canonical join representation of an element considered within a lower interval of a lattice is the same as its canonical join representation on the whole lattice.
Thus, we can obtain the canonical join representation of $T$ in  $\T_\skew$ by finding the canonical join representation of $\pi_T$ in  $[\id,\sigma_{\skew}]$, which is also the canonical join representation of $\pi_T$ in the weak order. By Theorem~\ref{thm:Reading_cjr_perm},
this is 
$$\pi_T=\bigvee \{\oldlambda(\pi_T,i): i\in\Des(\pi_T)\}.$$
Using that $\oldlambda(\pi_T,i)=\pi_{\La^{-1}(\{P_i\})}$ by Lemma~\ref{lem:oldlambdatab}, and $\Des(\pi_T)=\Asc(T)$ by Proposition~\ref{prop:interval}, we can rewrite this canonical join representation as 
$$\pi_T=\bigvee \{\pi_{\La^{-1}(\{P_i\})}: i\in\Asc(T)\}
=\bigvee \{\pi_{\La^{-1}(\{P\})}: P\in\La(T)\}.
$$
Finally, we use the isomorphism $\pi\mapsto\pi_T$ to obtain the stated canonical representation of~$T$.
    \end{proof}

\begin{corollary} \label{cor:join_irr} The join-irreducible elements of $\T_{\skew}$ are the tableaux  $T\in\SYT(\skew)$ for which $\card{\La(T)}=1$; or equivalently, $\asc(T)=1$. 
\end{corollary}

\begin{proof}
In any semidistributive lattice $L$, an element $x\in L$ is join-irreducible if and only if it has exactly one canonical joinand~\cite[Prop.~2.1]{Reading2015}. By Lemma~\ref{lem:cjr_tab}, the number of canonical joinands of $T$ equals the cardinality of $\La(T)$. Thus, $\asc(T)=\card{\La(T)}=1$ if and only if $T$ is join-irreducible in $\T_\skew$. 
\end{proof}

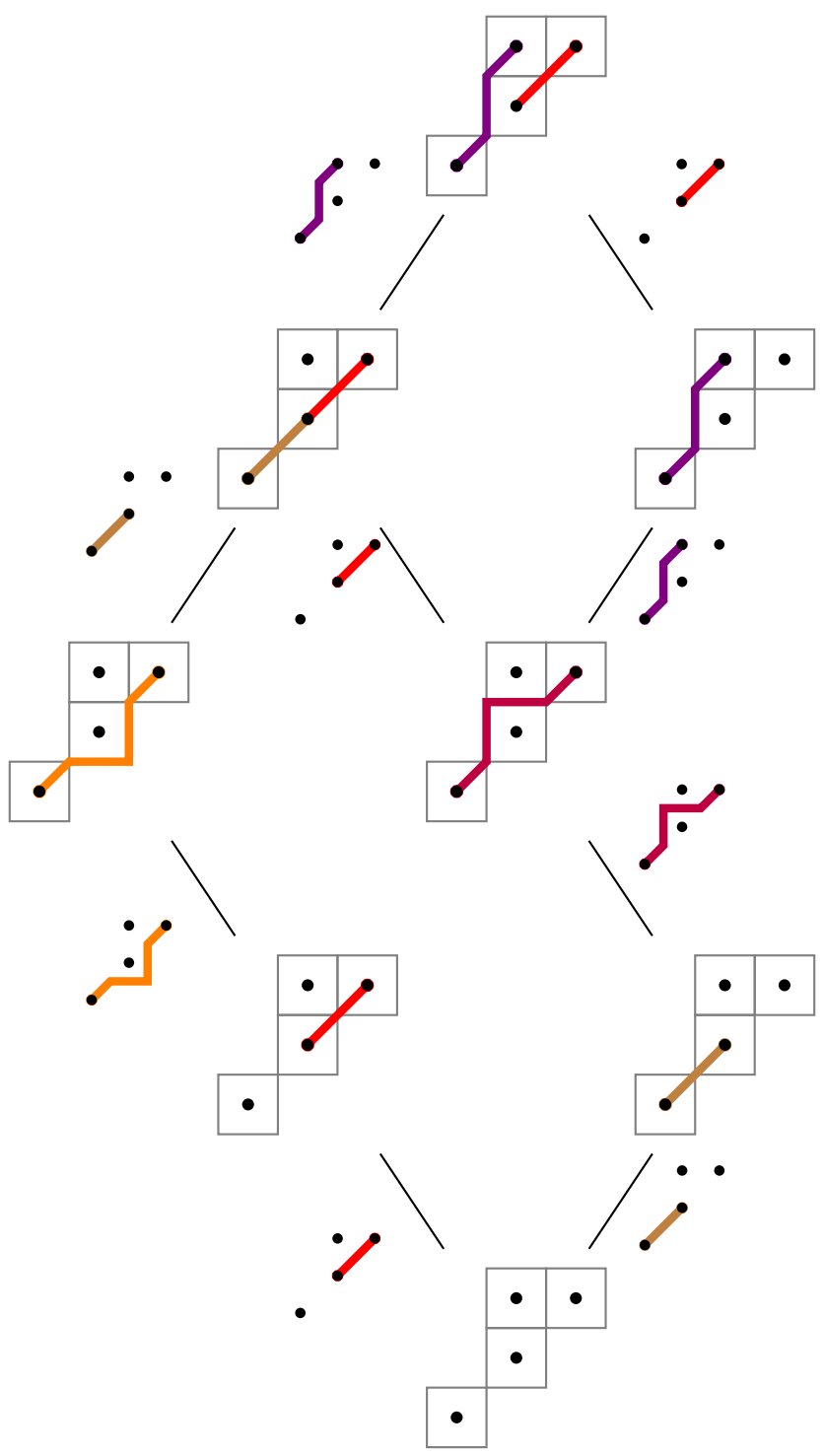
\begin{figure}[hbt]
\centering
\begin{tikzpicture}[scale=1.4]

\newcommand\tableau[1]{
\begin{tikzpicture}[scale=.4]
    \drawskewtab{3}{{1,1,0}}{{{,},{ },{ }}}
    #1
    \foreach \x in {(1,1),(2,2),(2,3),(3,3)}{
    \fill \x+(-.5,-.5) circle (.1);
    }
\end{tikzpicture}
}
\newcommand\edgelabel[1]{
\begin{tikzpicture}[scale=.25]
    #1
      \foreach \x in {(1,1),(2,2),(2,3),(3,3)}{
    \fill \x+(-.5,-.5) circle (.14);
    }
\end{tikzpicture}
}

\node (1) at (1,0) {\tableau{}};
\node (2) at (0,1.5) {\tableau{\drawpath{red}{(1.5,1.5)}{}}};
\node (3) at (2,1.5) {\tableau{\drawpath{brown}{(.5,.5)}{}}};
\node (4) at (-1,3) {\tableau{\drawpath{orange}{(.5,.5)}{\e\n}}};
\node (5) at (1,3) {\tableau{\drawpath{purple}{(.5,.5)}{\n\e}}};
\node (6) at (0,4.5) 
{\tableau{\drawpath{brown}{(.5,.5)}{}\drawpathsmall{red}{(1.5,1.5)}{}}};
\node (7) at (2,4.5) 
{\tableau{\drawpath{violet}{(.5,.5)}{\n}}};
\node (8) at (1,6) {\tableau{\drawpath{violet}{(.5,.5)}{\n}\drawpathsmall{red}{(1.5,1.5)}{}}};

\draw (1) -- node[below left] {\edgelabel{\drawpath{red}{(1.5,1.5)}{}}} (2);
\draw (1) -- node[right] {\edgelabel{\drawpath{brown}{(.5,.5)}{}}} (3);
\draw (2) -- node[below left] {\edgelabel{\drawpath{orange}{(.5,.5)}{\e\n}}} (4);
\draw (3) -- node[above right] {\edgelabel{\drawpath{purple}{(.5,.5)}{\n\e}}} (5);
\draw (4) -- node[above left] {\edgelabel{\drawpath{brown}{(.5,.5)}{}}} (6);
\draw (5) -- node[left] {\edgelabel{\drawpath{red}{(1.5,1.5)}{}}} (6);
\draw (5) -- node[right] {\edgelabel{\drawpath{violet}{(.5,.5)}{\n}}} (7);
\draw (6) -- node[above left] {\edgelabel{\drawpath{violet}{(.5,.5)}{\n}}} (8);
\draw (7) -- node[above right] {\edgelabel{\drawpath{red}{(1.5,1.5)}{}}} (8);

\end{tikzpicture}
\caption{The lattice from Figure~\ref{fig:T321_11} where each element $T$ is encoded by the family of paths $\La(T)$, and each edge is labeled by a join-irreducible as described in Definition~\ref{def:labels}. }
\label{fig:T321_11_paths}
\end{figure}

By Corollary~\ref{cor:join_irr}, join-irreducible elements in $\T_\skew$ are tableaux with exactly one ladder.
Thus, we can identify each join-irreducible element with its unique ladder, which is a path in $\P_\skew$.
In particular, the labels of the edges of $\T_\skew$ given by Definition~\ref{def:labels} can be viewed as paths in $\P_\skew$. 
Figure~\ref{fig:T321_11_paths} shows the edge labels of the lattice in Figure~\ref{fig:T321_11}.

\begin{lemma}
\label{lem:uc}
    For any $T\in\SYT(\skew)$, we have $\uplabels(T)=\{\La^{-1}(\{P\}):P\in\Ch(T)\}$.
\end{lemma}

\begin{proof}
For fixed $T$, covering relations of the form $T\lessdot U$ are in one-to-one correspondence with elements $i\in\HDes(T)$, that is, entries $i$ such that $i$ is northeast of $i+1$ in $T$. Because covering relations are defined by Bender--Knuth involutions, $T$ and $U$ differ only by swapping $i$ and $i+1$, and so $i\in\Asc(U)$. Moreover, by Definition \ref{def:chutes&ladders}, the path $P_i\in\Ch(T)$ from the cell containing $i$ to the cell containing $i+1$ whose north and east steps are determined by the boundary of $T^{<i}$ is the same as the path $P_i\in\La(U)$.

By Definition~\ref{def:labels}, the edges of the Hasse diagram between $U$ and the elements covered by $U$ are labeled with the canonical joinands of $U$.  
By Lemma~\ref{lem:cjr_tab}, these canonical joinands correspond to paths in $\La(U)$.
By the canonical join representation of the weak order (Theorem~\ref{thm:Reading_cjr_perm}), we have $\pi_T\vee \oldlambda(\pi_U,i)=\pi_{U}$. By Lemma~\ref{lem:oldlambdatab}, we have $\oldlambda(\pi_U,i)=\pi_{\La^{-1}(\{P_i\})}$, and so $\pi_T\vee\pi_{\La^{-1}(\{P_i\})}=\pi_{U}$. By Proposition~\ref{prop:interval}, it follows that $T\vee \La^{-1}(\{P_i\}) = U$. We deduce from Definition~\ref{def:labels} that the label $\ell(TU)$ is the join-irreducible tableau $\La^{-1}(\{P_i\})$.
\end{proof}

\begin{proof}[Proof of Theorem~\ref{thm:kappa}]
Given $T\in\SYT(\skew)$, we wish to show that $\rho(T)=\kappa(T)$. 
Using that $\rho=\Ch^{-1}\circ\La$ from Definition~\ref{def:row}, and the fact that $\Ch$ is a bijection by Corollary~\ref{cor:Ch}, this is equivalent to showing that $\La(T)=\Ch(\kappa(T))$.

Lemma~\ref{lem:cjr_tab}, together with Definition~\ref{def:labels}, implies that $\lowlabels(T)=\{\La^{-1}(\{P\}):P\in\La(T)\}$. By Definition~\ref{def:kappa}, we have $\lowlabels(T)=\uplabels(\kappa(T))$. 
And by Lemma~\ref{lem:uc}, we have $\uplabels(\kappa(T))=\{\La^{-1}(\{P\}):P\in\Ch(\kappa(T))\}$. We deduce that $\Ch(\kappa(T))=\La(T)$. 
\end{proof}

One may check Theorem~\ref{thm:kappa} in our running example by computing $\rho$ on the lattice in Figure~\ref{fig:T321_11} and comparing to the computation of $\kappa$ at the end of Section~\ref{sec:kapparow}.

\subsection{Consequences of the main theorem}

Rowmotion on the weak order on $\S_n$ is defined in \cite{Barnard2019} using the general construction from Definition~\ref{def:kappa}, and denoted by $\kappa$.
As a consequence of Theorem~\ref{thm:kappa}, we will prove in Corollary~\ref{cor:rubey} that this operation has an alternative description in terms of Reading's noncrossing arc diagrams.
 We start by precisely relating arc diagrams of permutations and arc diagrams determined by the ladders of tableaux in $\SYT(\delta_n/\delta_{n-1})$.
Recall from Corollary~\ref{cor:disjoint} that the map $T\mapsto\pi_T$ is a bijection between $\SYT(\delta_n/\delta_{n-1})$ and~$\S_n$.

\begin{lemma}\label{lem:arc_diagrams_permutations}
    For $T\in\SYT(\delta_n/\delta_{n-1})$, the noncrossing arc diagram of shape $\delta_n/\delta_{n-1}$ determined by the arcs in $\La(T)$ (as defined in Remark~\ref{rem:arc_diagrams})
    is equivalent to the noncrossing arc diagram of the permutation $\pi_T$ (as defined in Section~\ref{sec:arc_diagrams}), via the equivalence described in Remark~\ref{rem:arc_diagrams-connection} (i.e., $45^\circ$ rotation followed by vertical reflection).
\end{lemma}    

\begin{proof}
    
    We argue that for each $i\in\Asc(T)=\Des(\pi_T)$, the arc in the noncrossing arc diagram of shape $\delta_n/\delta_{n-1}$ determined by the path $P_i\in\La(T)$ is equivalent to the arc in the noncrossing arc diagram for the permutation $\pi_T$ determined by the descent $\pi_T(i)>\pi_T(i+1)$. 
    
    Indeed, by Definition~\ref{def:chutes&ladders}, the path $P_i\in\La(T)$ goes from the cell containing $i$ to the cell containing $i+1$ in $T$, and leaves to its left the cells in $T^{<i}$. Suppose that the labels of the cells containing $i$ and $i+1$, in the labeling by rows of $\delta_n/\delta_{n-1}$, are $b=\pi_T(i)$ and $a=\pi_T(i+1)$. Then the arc diagram of $\pi_T$ has an arc from $b$ to $a$ that leaves to its left the dots with labels $k\in(a,b)$ such that the entry $k$ appears to the left of $a$ and $b$ in $\pi_T$. But $k$ appears to the left of $a$ and $b$ in $\pi_T$ precisely if it is the label of a cell in $T^{<i}$.
\end{proof}

\begin{figure}[htbp]
\begin{center}
\begin{tikzpicture}[scale=.5]

\begin{scope}[shift={(4,10)},
    dot/.style={circle,fill=black,inner sep=1.7pt}]

    \draw[red,thick] (0,2) to [bend left=55] (0,5);
    
    \draw[red,thick] (0,3)  to [bend right=55] (0,5.5) to [bend left=55] (0,8);

  \foreach \i in {1,...,8} {
    \node[dot] at (0,\i) {};
    \node[left=3pt] at (0,\i) {\(\i\)};
  }
    \draw (2,-.5) node[left] {$\pi=\pi_T=15248367$};
    \draw (2,-1.5) node[blue,left] {$\pi_T^{-1}=13642785$};
\end{scope}

\begin{scope}[shift={(15,10)},
    dot/.style={circle,fill=black,inner sep=1.7pt}]

    \draw[red,thick] (0,2) to [bend right=55] (0,5);
    
    \draw[red,thick] (0,3)  to [bend left=55] (0,5.5) to [bend right=55] (0,8);

  \foreach \i in {1,...,8} {
    \node[dot] at (0,\i) {};
    \node[left=3pt] at (0,\i) {\(\i\)};
  }
    \draw (2,-.5) node[left]  {$\tau=\pi_U=16783452$};
    \draw (2,-1.5) node[blue,left] {$\pi_U^{-1}=18567234$};
\end{scope}

\begin{scope}[shift={(26,10)}]
    \draw (2,-.5) node[left] {$\tau^{\rev}=\pi_V=25438761$};
    \draw (2,-1.5) node[blue,left] {$\pi_V^{-1}=81432765$};
\end{scope}

\drawcircpath{red}{(3.5,3.5)}{\n\e\n\e}
\drawcircpath{violet}{(.5,.5)}{\n\e\n\e\e\n\e\n}
\drawskewtab{8}{{7,6,5,4,3,2,1,0}}{{{1},{3},{6},{4},{2},{7},{8},{5}}}
\draw (2,4) node[left]{$T=$};
\draw[->] (8.5,4) -- node[above]{$\flip$} (9.5,4);
\begin{scope}[shift={(11,0)}]
\drawcircpath{red}{(3.5,3.5)}{\e\n\e\n}
\drawcircpath{violet}{(.5,.5)}{\e\n\e\n\n\e\n\e}
\drawskewtab{8}{{7,6,5,4,3,2,1,0}}{{{1},{8},{5},{6},{7},{2},{3},{4}}}
\draw (2,4) node[left]{$U=$};
\draw[->] (8.5,4) -- node[above]{$\comp$} (9.5,4);
\end{scope}
\begin{scope}[shift={(22,0)}]
\drawcircpath{red}{(3.5,3.5)}{\n\e\n\e}
\drawcircpath{violet}{(.5,.5)}{\n\e\n\e\e\n\e\n}
\drawskewtab{8}{{7,6,5,4,3,2,1,0}}{{{8},{1},{4},{3},{2},{7},{6},{5}}}
\draw (2,4) node[left]{$V=$};
\end{scope}
\end{tikzpicture}
\end{center}
\caption{Example of the computation of rowmotion on $\pi=\pi_T\in\S_n$ by applying $\comp\circ\flip$ to $T\in\SYT(\delta_n/\delta_{n-1})$. The paths shown on $T$ and $U$ tableaux are ladders, while those on $V$ are chutes.}
\label{fig:weak_order_row}
\end{figure}
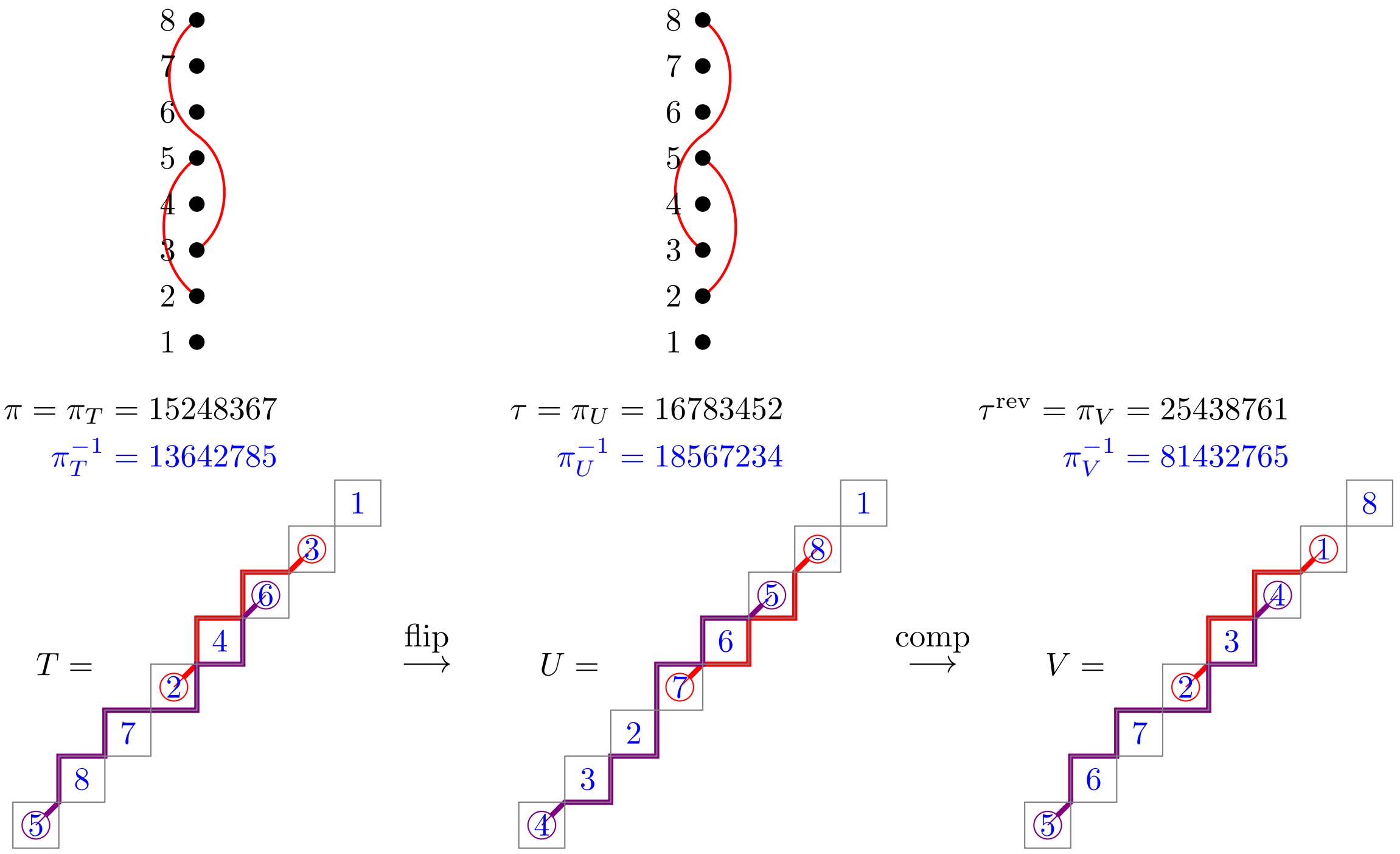

The {\em reversal} of a permutation $\pi=\pi(1)\pi(2)\dots\pi(n)$ is the permutation $\pi^{\rev}=\pi(n)\dots\pi(2)\pi(1)$.
The following alternative description of rowmotion on the weak order was conjectured by Rubey~\cite{FindStatMp00326}. Our proof is illustrated with an example in
Figure~\ref{fig:weak_order_row}.

\begin{corollary}\label{cor:rubey}
    Rowmotion on the weak order can be computed as follows: given $\pi\in\S_n$, reflect the arc diagram of $\pi$ horizontally and apply reversal to the resulting permutation to obtain $\kappa(\pi)$.
\end{corollary}

\begin{proof}
Let $\pi\in\S_n$, and let $\tau$ be the permutation whose arc diagram is obtained by horizontally reflecting the arc diagram of $\pi$. We want to show that $\kappa(\pi)=\tau^{\rev}$.

By Corollary~\ref{cor:disjoint}, the map $T\mapsto \pi_T$ is a poset isomorphism between $\T_{\delta_n/\delta_{n-1}}$ and the weak order on $\S_n$. Thus, we can assume that $\pi=\pi_T$ for some $T\in\SYT(\delta_n/\delta_{n-1})$, compute rowmotion on $T$, and use the fact that $\kappa(\pi)=\pi_{\kappa(T)}$. By Theorem~\ref{thm:kappa}, $\kappa(T)=\rho(T)$, so we can use the description of rowmotion as $\row=\Ch^{-1}\circ\La$ from Definition~\ref{def:row}.

By Lemma~\ref{lem:arc_diagrams_permutations}, the arc diagram of shape $\delta_n/\delta_{n-1}$ determined by the arcs in $\La(T)$ is equivalent to the arc diagram of $\pi_T$. Because of this equivalence, horizontally reflecting the arcs in the arc diagram of $\pi_{T}$
corresponds to reflecting the family of paths $\La(T)$ across the northeast-to-southwest diagonal.
Let $\flip:\SYT(\delta_n/\delta_{n-1})\rightarrow\SYT(\delta_n/\delta_{n-1})$ denote the map that reflects $\La(T)$ across this diagonal and then uses $\La^{-1}$ (as in Definition~\ref{def:La-inv}) to produce the tableau $U=\flip(T)$.

Define another map $\comp:\SYT(\delta_n/\delta_{n-1})\rightarrow\SYT(\delta_n/\delta_{n-1})$ by letting $\comp(U)=V$, where $V$ is the unique tableau such that $\pi_V$ is the reversal of $\pi_U$.
This property is equivalent to $\pi_U(k)=\pi_V(n+1-k)$, which says that the cell containing $k$ in $U$ is the cell containing $n+1-k$ in $V$, for all $k$. Thus, the map $\comp$ simply replaces each tableau entry $k$ with its complementary value $n+1-k$.

Let $V=\comp(\flip(T))$. We wish to show that $\row(T)=V$, or equivalently, that 
$\La(T)=\Ch(V)$, that is, the ladders of $T$ equal the chutes of $V$. 

Let $P\in\La(T)$, and let $C$ and $C'$ be the cells at its endpoints. Let $i$ be the entry in cell $C$ so that $i+1$ is the entry in cell $C'$. Then the cells to the left of $P$ have entries $<i$, and the cells to the right of $P$ have entries $>i+1$.

By construction, the tableau $U=\flip(T)$ has a ladder $\widehat{P}$ from $C$ to $C'$ that is obtained by reflecting $P$. Thus, $U$ still has an ascent between these two cells. Let $j$ be the entry in cell $C$ of $U$, so that $j+1$ is the entry in cell $C'$. Then the cells to the left of $P$ (which are those to the right of $\widehat{P}$) have entries $>j+1$, and the cells to the right of $P$ (which are those to the left of $\widehat{P}$) have entries $<j$.

The tableau $V=\comp(U)=\comp(\flip(T))$ is obtained from $U$ by replacing each entry $k$ with $n+1-k$. Thus, $C$ contains entry $n+1-j$ and $C'$ contains entry $n-j$. Additionally, the cells to the left of $P$ have entries $<n-j$, and the cells to the right of $P$ have entries $>n+1-j$. This means that $P\in\Ch(V)$.

We have shown that every ladder of $T$ is a chute of $V$. A similar argument (applying $\flip\circ\comp$ to $V$) shows that every chute of $V$ is a ladder of $T$. Thus, $\La(T)=\Ch(V)$.
\end{proof}

We end this subsection showing that, for some shapes $\skew$, rowmotion on $\SYT(\skew)$ has the same orbit structure as promotion.

\begin{corollary}
\label{cor:match}
Suppose that $\skew$ has at most two rows or at most two columns, or is contained in a hook or a reverse hook. Then rowmotion and promotion on $\SYT(\skew)$ have the same orbit structure. 
\end{corollary}
\begin{proof}
By Theorem~\ref{thm:kappa}, rowmotion on $\SYT(\skew)$ can be described using either $\rho$ (from Definition~\ref{def:row}) or $\kappa$ (from Definition~\ref{def:kappa}), since they are the same map. For this proof, we use its description as $\kappa$.

For the shapes $\skew$ considered in the statement, $\T_{\skew}$ is a distributive lattice by Theorem~\ref{thm:dist}, and hence it is isomorphic to the lattice of order ideals of its subposet of join-irreducible elements. By Remark~\ref{remark:rowequalsrow}, the rowmotion map $\kappa$ on $\T_{\skew}$ is equivalent to order ideal rowmotion on this subposet.

Suppose first that $\skew=(m^2)$. The subposet of join-irreducibles of $\T_{(m^2)}$ is isomorphic to $\A_{m-1}$ (defined in Section~\ref{sec:two_row}), as shown in \cite[Sec.\ 4.3]{prorow} and \cite[Sec.\ 3.5]{EC1}. 
Rowmotion on the distributive lattice $\T_{(m^2)}$ corresponds to order ideal rowmotion on $\A_{m-1}$, which can be computed by toggling from top to bottom. Promotion on $\SYT(m^2)$ corresponds to toggling order ideals of $\A_{m-1}$ by columns from left to right \cite[Thm.~4.12]{prorow}.  By \cite[Thm.~5.2]{prorow}, the elements of the toggle group obtained by toggling from top to bottom and from left to right are conjugate of each other, and thus they have the same orbit structure. 

For skew shapes $\skew$ contained in $(m^2)$, the subposet of join irreducibles of $\T_\skew$ is isomorphic to a subposet of $\A_{m-1}$, by \cite[Thms.\ 4.12 and 5.2]{prorow}, and a similar argument applies.
    
If $\skew$ has at most $2$ columns, then the posets $\T_{\skew}$ and $\T_{\lambda'/\mu'}$ are dual to each other by Lemma~\ref{lem:dual}. Rowmotion on the dual of a semidistributive lattice is the inverse of rowmotion on the original lattice. Thus, rowmotion has the same orbit structure on $\T_{\skew}$ and $\T_{\lambda'/\mu'}$. The shape $\lambda'/\mu'$ has at most two rows, so it is covered in the previous case.  

Finally, if $\skew$ is contained in a hook or reverse hook then $\T_{\skew}$ is isomorphic to $\T_{(a+b,b) /(b)}$ for some $a$ and $b$, as explained in the proof of Theorem~\ref{thm:dist}, which is also covered in the previous case.
\end{proof}

\section{Open problems and future work}\label{sec:open}

\subsection{Orbits of rowmotion}
\label{sec:orbits}
It is natural to ask for what shapes $\skew$ is the order of rowmotion on $\SYT(\skew)$ relatively small compared to $n=|\skew|$. In the case $\skew=\delta_n/\delta_{n-1}$, which corresponds to the weak order on $\S_n$, it has been observed \cite[Sec.\ 7.11]{ThomasWilliams2019} that the orbit sizes and the order of rowmotion grow as $n$ increases. 
Table~\ref{tab:Rowmotion_order} suggests that staircases and rectangles display a similarly wild behavior.

\renewcommand{\arraystretch}{1}

\begin{table}[htbp]
    \centering
    \begin{tabular}{c|r|r}
       $\lambda$  &  Largest orbit size & Order \\
       \hline\hline
        (2,1) & 2 & 2 \\
               \hline
        (3,2,1) & 6 & 12\\
               \hline
        (4,3,2,1) & 100 & 25,225,200\\
        \hline
        (3,3,3) & 8 & 8\\
        \hline
        (4,4,4) & 66 & 112,200
    \end{tabular}
    \caption{Data on the orbits of rowmotion on $\SYT(\lambda)$ for some staircase and rectangular shapes.}
\label{tab:Rowmotion_order}
\end{table}

On the other hand, for rowmotion on the weak order, it is shown in \cite[Sec.\ 6.2]{HopkinsCDE} that some statistics are {\em homomesic}, meaning that they have the same average on each orbit. One may ask whether this phenomenon extends to rowmotion on $\SYT(\skew)$ for other skew shapes $\skew$.

The converse of Corollary~\ref{cor:match} is not true, as there are other shapes $\skew$ for which promotion and rowmotion on $\SYT(\skew)$ have the same orbit structure. The smallest such shapes have four cells. For example, for the skew shape of Figure~\ref{fig:T321_11}, both promotion and rowmotion have one orbit of size $8$. On the other hand, as shown in Table~\ref{tab:Rowmotion_order}, rowmotion has order $8$ for $\lambda=(3,3,3)$, while promotion is known to have order $9$. Similarly, for $\lambda=(4,4,4)$, rowmotion has order $112,200$, while promotion has order $12$; for $\lambda=(4,3,2,1)$, rowmotion has order $25,225,200$, while promotion has order $20$. It would be interesting to 
    determine all shapes for which the orbit structure of rowmotion on $\SYT(\skew)$ matches that of promotion.

\subsection{Generalization to other posets}

In this section we discuss possible generalizations of our results about standard Young tableaux  to linear extensions of other posets.

Let $P$ be any poset with $n$ elements. A {\em linear extension} of $P$ is an order-preserving bijection $\tau:P\to\{1,2,\dots,n\}$.
Fix a particular linear extension $\omega$, called a {\em natural labeling}, and identify each element $x$ of $P$ with its label $\omega(x)$. Every linear extension $\tau$ can then be identified with the permutation obtained by reading the labels of the poset in the order given by the linear extension, that is, the permutation $\omega(\tau^{-1}(1))\dots\omega(\tau^{-1}(n))\in\S_n$. We denote the set of these permutations by $\L(P,\omega)$. 
The descents of a linear extension $\tau$ with respect to $\omega$ are defined as the descents of the associated permutation, that is,
\[
\Des_\omega(\tau):=\{i:\,\omega(\tau^{-1}(i))>\omega(\tau^{-1}(i+1))\}.
\]
It is known that the distribution of $\Des_\omega$ does not depend on the natural labeling~$\omega$.

This paper has focused on the case of posets that arise from skew shapes $\skew$. Specifically, let $P_\skew$ be the poset whose elements are $\{(r,c):1\le r\le k,\, \mu_r<c\le\lambda_r\}$ ordered coordinate-wise. Then, linear extensions of $P_\skew$ correspond to standard Young tableaux in $\SYT(\skew)$. By fixing the natural labeling by rows as in Definition~\ref{def:labeling_by_rows}, each tableau $T\in\SYT(\skew)$ is identified with the permutation $T_R(T^{-1}(1))\dots T_R(T^{-1}(n))$, which we denoted by $\pi_T$ in Section~\ref{sec:SYTposet}. Proposition~\ref{prop:interval} showed that the resulting set of permutations $\L(P_\skew,T_R)$ forms an interval in the weak order on $\S_n$ and that descents of such permutations correspond to ascents of the associated tableaux.

For an arbitrary poset $P$, we can order the elements of $\L(P,\omega)$ by defining $\pi$ to cover $\widehat\pi$ if $\widehat\pi$ is obtained from $\pi$ by swapping two adjacent entries that form a descent of $\pi$. At the level of the linear extension $\tau$, a covering relation corresponds to swapping two consecutive values $i$ and $i+1$ associated to incomparable elements. In general, the resulting poset on $\L(P,\omega)$ may not be an interval in the weak order. In fact, it may have multiple maximal elements, and it may not be a lattice.

Bj\"orner and Wachs showed in \cite[Sec.~6]{Bjorner_Wachs-permutations} that posets $P$ that have a labeling $\omega$ such that the resulting poset on $\L(P,\omega)$ is an interval in the weak order are precisely two-dimensional posets (that is, those that can be described as the intersection of two linear orders).
Additionally, they characterize in \cite[Thm.~6.8]{Bjorner_Wachs-permutations} which labelings $\omega$ have this property, and call them {\em regular labelings}. A natural labeling $\omega$ is regular if for all $x,y,z\in P$ such that $x<z$ and $\omega(x)<\omega(y)<\omega(z)$, we have $x<y$ or $y<z$.

It follows that, given any two-dimensional poset $P$ and any regular natural labeling $\omega$, the resulting poset on $\L(P,\omega)$ is a lower interval in the weak order on $\S_n$. In particular, it is a semidistributive lattice, so one can define a rowmotion operation $\kappa$ as in Definition~\ref{def:kappa}. This suggests the following problem.

\begin{problem}
    Generalize our definition of chutes and ladders on standard Young tableaux to other two-dimensional posets, and use it to give a direct description of rowmotion on linear extensions of these posets, similar to our construction in Section~\ref{sec:rowmotion-chutes&ladders}.
\end{problem}

A class of two-dimensional posets of particular interest are forests, which are posets where each element is covered by at most one element. Forests play an important role in~\cite{Bjorner_Wachs,Bjorner_Wachs-permutations}, where it is shown that postorder is a regular labeling. 
Recall that postorder is the ordering of the nodes of a tree obtained by recursively reading the subtrees of the root from left to right, followed by the root itself.

We conclude with some questions about equidistributed statistics on linear extensions. For a given natural labeling $\omega$ of a poset $P$, consider the statistic $\des_\omega(\tau)=\card{\Des_\omega(\tau)}$ on linear extensions $\tau$ of $P$. 
It is shown in \cite[Thm.~3.15.8]{EC1}, using the theory of $P$-partitions, that the distribution of $\des_\omega$ on linear extensions of $P$ does not depend on the choice of $\omega$. More generally, it can be shown  using Gessel's quasisymmetric functions~\cite{gessel_multipartite_1984} that the same holds for the set-valued statistic $\Des_\omega$. A recursively-defined bijection taking $\Des_\omega$ to $\Des_{\omega'}$, for any two natural labelings $\omega$ and $\omega'$, is given in \cite[Prop.~7.2]{EliSYT}. However, no direct bijective proof of the fact that these statistics are equidistributed is known in general. 

For posets of the form $P_\skew$, one can view the statistic $\Des_\omega$ as a generalization of both $\Asc$ and $\HDes$. Indeed, if $\omega$ is the natural labeling by rows and $\omega'$ is the natural labeling by columns (see Definition~\ref{def:labeling_by_rows}), then $\Des_\omega=\Asc$ and $\Des_{\omega'}=\HDes$. Thus, our rowmotion operation $\rho$ on $\SYT(\skew)$ gives a direct bijective proof of the equidistribution of the statistics $\des_\omega=\asc$ and $\des_{\omega'}=\hdes$. One can try to generalize this to other settings.

\begin{problem}
For a given poset $P$ and two natural labelings $\omega$ and $\omega'$, find a direct bijection that takes $\des_\omega$ to $\des_{\omega'}$.
\end{problem}

The definitions of $\AscC$ and $\HDesC$ from Section~\ref{sec:chutes&ladders} can also be generalized to other posets, by defining $$\DesC_\omega(\tau):=\{\{\tau^{-1}(i),\tau^{-1}(i+1)\}:\,\omega(\tau^{-1}(i))>\omega(\tau^{-1}(i+1))\}$$
for any natural labeling $\omega$ of a poset $P$ and any linear extension $\tau$. Intuitively, this statistic records the pairs of elements in $P$ where the descents of the linear extension occur.
However, unlike in the case of $\Des_\omega$, the distribution of the statistic $\DesC_\omega$ over linear extensions of $P$ does depend on the natural labeling $\omega$ in general.

\subsection*{Acknowledgements}
The first author was partially supported by Simons Collaboration Grant \#929653. The second author was partially supported by NSF Grants DMS-2247089 and DMS-2554137 and Simons Foundation Gift MP-TSM-00002802. The authors thank Ben Adenbaum, Emily Barnard, Sam Hopkins, Martin Rubey, and Nathan Williams for helpful conversations.
SageMath was helpful for some computations. No AI has been used in the creation of this paper, other than to make some figures and to check for typos. 

\bibliographystyle{plain}
\bibliography{rowmotionSYT}

\end{document}